\documentclass[a4paper,12pt]{article}
\usepackage{amsfonts,amsmath,amssymb}
\usepackage{graphicx,subfigure}

\usepackage{epsfig}
\usepackage{verbatim}
\usepackage{color}
\usepackage{booktabs}
\usepackage{siunitx}
\usepackage[geometry]{ifsym}
\usepackage{float}
\usepackage{caption}
\usepackage{blindtext}
\usepackage{fancyhdr}

\usepackage{url}
\urldef{\mailsa}\path|{cesare.bracco, carlotta.giannelli, alessandra.sestini}@unifi.it|
\def \d {{\bf d}}

\def\RR{{\mathop{{\rm I}\kern-.2em{\rm R}}\nolimits}}
\def\PP{{\mathop{{\rm I}\kern-.2em{\rm P}}\nolimits}}

\def\NN{{\mathop{{\rm I}\kern-.2em{\rm N}}\nolimits}}

\newcommand{\be}{\begin{equation}}
\newcommand{\ee}{\end{equation}}
\newcommand{\ba}{\begin{eqnarray}}
\newcommand{\ea}{\end{eqnarray}}
\newcommand{\bi}{\begin{itemize}}
\newcommand{\ei}{\end{itemize}}

\newcommand{\supp}{\mathop{\mathrm{supp}}}

\newcommand{\myspan}{\mathop{\mathrm{span}}}
\newcommand{\trunc}{\mathop{\mathrm{trunc}}}
\newcommand{\Trunc}{\mathop{\mathrm{Trunc}}}

\newcommand{\diam}{\mathop{\mathrm{diam}}}

\newtheorem{dfn}{Definition}
\newtheorem{remark}{Remark}

\newtheorem{exm}{Example}

\usepackage{array}

\newcolumntype{C}[1]{>{\centering\let\newline\\\arraybackslash\hspace{0pt}}m{#1}}

\def\M3AS{Math.\ Models\ Methods\ Appl.\ Sci.}

\begin{document}

%\begin{frontmatter}

\title{Adaptive scattered data fitting by extension of local approximations to hierarchical splines}

\author{Cesare Bracco\thanks{cesare.bracco@unifi.it}, Carlotta Giannelli\thanks{carlotta.giannelli@unifi.it}, Alessandra Sestini \thanks{alessandra.sestini@unifi.it}\\
\small Dipartimento di Matematica e Informatica ``U. Dini'', Universit\`a degli Studi di Firenze, \\
\small Viale Morgagni 67/a, 50134 Firenze,  Italy \\
%\small \texttt{cesare.bracco@unifi.it, carlotta.giannelli@unifi.it, alessandra.sestini@unifi.it}\\
%\and
}
\date{}
\maketitle

\begin{abstract}
We introduce an adaptive scattered data fitting scheme as extension of local least squares approximations to hierarchical spline spaces. To efficiently deal with non-trivial data configurations, the local solutions are described in terms of (variable degree) polynomial approximations according not only to the number of data points locally available, but also to the smallest singular value of the local collocation matrices. These local approximations are subsequently combined without the need of additional computations with the construction of hierarchical quasi-interpolants described in terms of truncated hierarchical B-splines. A selection of numerical experiments shows the effectivity of our approach for the approximation of real scattered data sets describing different terrain configurations.
\end{abstract}

{\bf Keywords:} Scattered data fitting, Hierarchical splines, THB-splines, Local least squares, Quasi-Interpolation. 
%\end{frontmatter}

\section{Introduction}

Surface reconstruction of unstructured large data sets requires suitable \emph{adaptive} schemes that facilitate the computation of high-quality approximations with an increased level of resolution only in strictly localized areas. The resulting compact representation automatically identifies the parts of the domain where an increased number of degrees of freedom is needed according to the data distribution (high concentrations of data points usually define local details to be suitably reconstructed). The problem of reconstructing scattered data of high complexity arises in various application areas, ranging from scanner acquisitions to geographic benchmarks, and it is a relevant component for industrial and medical purposes where the visualization and subsequent manipulation of large random data configurations are usually required. Note that scattered data can have significantly different distributions, e.g. data with highly varying density, data with voids, contour data.

Standard computer aided design software tools based on the tensor-product B-spline model do not provide local refinement capabilities. Spline adaptivity may easily be achieved by considering multilevel B-spline extensions, where the tensor-product structure is preserved at any level. Hierarchical B-spline constructions of this kind were originally proposed by \cite{forsey1988} to define hierarchical spline surfaces in terms of a sequence of overlays. By considering truncated hierarchical B-splines (THB-splines) \cite{giannelli2012}, it is possible to define a strongly stable basis for hierarchical spline spaces that forms a convex partition of unity \cite{giannelli2014}. THB-splines also guarantee the so--called \emph{preservation of coefficients}: truncated basis functions preserve the coefficients of functions expressed in terms of B-splines of a certain hierarchical level, see again \cite{giannelli2014}. This allows us to directly extend any quasi-interpolation operator defined in the space of tensor-product B-splines to the hierarchical setting without the need of additional computation \cite{speleers2016}. Note that the hierarchical B-spline model may be applied in combination with uniform and non-uniform refinement, different degrees and smoothness, and related constructions easily extend to the general multivariate setting.

We present an adaptive scattered data fitting scheme by extension of local discrete least squares approximations to hierarchical spline spaces. To efficiently deal with non-trivial data configurations, the local solutions are described in terms of variable degree polynomial approximations, according not only to the number of data points locally available, but also to the smallest singular value of the local collocation matrices. Note that the possibility of using higher degree polynomials in each local approximation can be exploited only for sufficiently dense data sets. These local approximations are subsequently converted in local B--spline form and directly combined with the construction of hierarchical quasi-interpolants described in terms of the truncated basis.  The choice of defining local approximations in terms of polynomials instead of splines makes our approach more robust when data with highly varying density or voids are considered. Within the hierarchical framework, this local quasi-interpolation operator is combined with an adaptive strategy that at each refinement step identifies a set of basis functions marked for refinement. In order to do this, we developed a marking strategy to identify the subset of hierarchical basis functions to be substituted by finer ones in the spline hierarchy. This allows us to locally improve the accuracy of the approximation according to the scattered data distribution.

The structure of the paper is as follows. Section~\ref{sec:relatedwork} presents a brief overview of some related works dealing with scattered data fitting, while Section~\ref{sec:thb} recalls the definition of THB-splines and related properties. Section~\ref{sec:scatteredqi} introduces the developed least squares local data-dependent quasi-interpolation approach, while Section~\ref{sec:adap} introduces the strategy for the definition of our final adaptive hierarchical approximation. A selection of numerical experiments for the approximation of real data sets describing different terrain configurations is presented in Section~\ref{sec:exm}. Finally, Section~\ref{sec:end} concludes the paper.

%%%%%%%%%%%%%%%%%%%%%%%%%%%%%%%%%%%%%%%%%%%%%%%%%%%%%%%%%%%%%%%%%%%%

\section{Related works}\label{sec:relatedwork}

One possible choice for scattered data fitting relies in methods based on radial basis functions (RBFs) since the space naturally depends on the local distribution and density of the input data, see e.g. \cite{fasshauer2007,wendland2004}. Another classic approach to the problem is represented by partition of unity methods, where (suitably chosen) non-negative, compactly supported and linearly independent functions, which form a partition of unity, are combined with local approximants of the data, see, e.g., \cite{bracco2012,fasshauer2007}.
%\cite{ohtake2005} 
Among others, a method that combines adaptive partition of unity with least squares fitting based on radial basis functions was proposed in \cite{ohtake2006}, while \cite{liu2012} presented a scattered data quasi-interpolation scheme based on RBFs.

Focusing on splines, a two-stage method for scattered data fitting by direct extension of local polynomials to bivariate splines was presented in \cite{davydov2004} by considering local  discrete least squares polynomial approximations \cite{davydov2002}. In order to increase the adaptivity of the method, this approach was further developed in \cite{davydov2005,davydov2006} by considering hybrid local approximations in terms of polynomials and of linear combinations of radial basis functions in the first stage of the method. In both cases the final  approximation is a spline (represented in local Bernstein form) on a regular triangulation (of a suitable extension) of the domain. Scattered data fitting can be naturally addressed with splines on irregular triangulations. Some recent theoretical studies in this context are, for example, \cite{nurnberger2005}, where interpolation is proposed for exact data, and \cite{lai2009} where a domain decomposition approach is introduced. In addition to discrete least squares, minimal energy and penalized least squares are considered in the first stage of the method in \cite{lai2009}, which can also be listed among the two-stage approaches for scattered data fitting. Furthermore, a quasi-interpolation scheme based on irregular triangulation has also been recently investigated \cite{speleers2015b}. Multilevel least squares approximation of scattered data over binary triangulations was presented in \cite{hjelle2005}. 

Moving to tensor-product splines, a more standard choice in computer aided applications than splines on triangulations, a two-stage method based on extended B-splines with focus on curvilinear domains  has  been recently introduced in the literature \cite{davydov2014}. However, this work provides theoretical results for sufficiently dense data, and, consequently, does not cover real scattered data sets. 
Interpolation and approximation of scattered {3D} data with hierarchical tensor-product B-splines was addressed in \cite{greiner1997} by following the hierarchical approach originally proposed in \cite{forsey1988}. The method was based on a global least squares minimization with fairing which is also the basic approximation choice adopted in \cite{rabut2005,kiss2014b}. More precisely, in \cite{rabut2005} local tensor-product functions on suitable subdomains are used via repeated knot insertion, while in [15] the THB-spline model is exploited. Interpolation and least square approximation of gridded data with hierarchical splines was also proposed in \cite{forsey1995} by taking advantage of the local tensor-product structure of any overlay of the hierarchical spline surface \cite{forsey1988}. Note that in the multilevel approach to hierarchical B-splines followed by \cite{forsey1988,forsey1995,greiner1997} the use of this kind of global scheme on any refinement level is motivated by the natural assumption of a limited number of degrees of freedom on a single level. When a hierarchical B-spline basis is constructed over the whole domain instead, the solution of the approximation problem is directly computed in the corresponding hierarchical spline space \cite{giannelli2012,kraft1997}. An adaptive extension of local approximation to hierarchical splines, as the two-stage method here proposed, is then the suitable choice in this case.
 Partial results related to multilevel B-spline quasi-interpolation of unstructured data were presented in \cite{lee2005}. This refinement construction, however,  defines a level-by-level correction in the tensor-product B-spline space that necessarily requires suitable data distributions. 
Consequently, general scattered data configurations with varying density cannot be covered with this approach. Approximation of large terrain data sets with locally refined B-splines was recently addressed \cite{skytt2015}.

%%%%%%%%%%%%%%%%%%%%%
\section{THB-splines and hierarchical quasi--interpolation}\label{sec:thb}

In this section we briefly introduce hierarchical B-spline spaces and summarize the construction of their truncated basis, which can be conveniently used for the definition of effective hierarchical spline quasi-interpolation operators.

Let 
\[
V^0 \subset V^1 \subset \ldots \subset V^{M-1} 
\qquad\text{and}\qquad
\left\{
{\cal B}_\d^{0}, {\cal B}_\d^{1},\ldots, {\cal B}_\d^{M-1}
\right\}
\]
be a sequence  of nested $r$-variate  tensor-product spline spaces $V^\ell$ defined on a hyper-rectangle $\Omega \subset \RR^r$ and a sequence of corresponding B-spline bases ${\cal B}_\d^{\ell}$ of degree $\d :=(d_1,\ldots,d_r)$, for $\ell=0,\ldots,M-1.$ In the following we indicate as  ${\cal G}^{\ell}$ the tensor-product grid associated to $V^{\ell}$ and as $\Gamma_{\d}^{\ell}$ the set of multi--indices with $r$ components such that  
\[
{\cal B}_\d^{\ell} := \left\{
B_J^{\ell},\, J\in \Gamma_{\d}^{\ell}
\right\}.
\]
We also consider a nested sequence of closed domains 
\[
\Omega^{0}\supseteq \Omega^{1}\supseteq \ldots \supseteq \Omega^{M},
\]
where any $\Omega^\ell\subset\Omega$ is defined as a collection of cells belonging to ${\cal G}^{\ell},$ for $\ell=0,\ldots, M-1,$ and $\Omega^M=\emptyset$. We define the hierarchical mesh ${\cal G}_{\cal H}$ as the set of all the {\it active} (namely, no further refined) cells belonging to $\Omega^{\ell},$ for each {\it level} $\ell=0,\ldots,M-1.$

The hierarchical B-spline basis is defined by iteratively selecting basis function at different levels,  according to the following construction \cite{kraft1997,vuong2011}.

\begin{dfn}\label{dfn:hb}
The hierarchical B-spline basis ${\cal H}_\d({\cal G}_{\cal H})$ of degree $\d$ 
with respect to the mesh ${\cal G}_{\cal H}$ is defined as
\begin{equation*}
{\cal H}_\d({\cal G}_{\cal H})  :=  \left\{
B_J^{\ell} \in {\cal B}_\d^\ell \, : \, J \in A_\d^{\ell}, \,\, \ell=0,...,M-1  \right\},
\end{equation*} 
where  
\begin{equation*} \label{acitve}
A_\d^{\ell} \,:=\, \{ J \in \Gamma_{\d}^{\ell}\,:\,  
\supp B_J^{\ell}  \subseteq \Omega^\ell\,\,\wedge \,\, 
\supp B_J^{\ell} \not\subseteq \Omega^{\ell+1} \}
\end{equation*}
is the active set of multi-indices $A_\d^{\ell} \subseteq  \Gamma_{\d}^\ell$ and $\supp B_J^{\ell}$ denotes the intersection of the support of $B_J^{\ell}$ with $\Omega^0$.
\end{dfn}

Hierarchical B-splines form a basis for the multilevel space $S_{\cal H} := \myspan {\cal H}_\d({\cal G}_{\cal H}) $ associated to the hierarchical mesh ${\cal G}_{\cal H}$.
 
Let $s \in V^\ell$ be represented in terms of B--splines of the refined space $V^{\ell+1}$ as
$$ s  = \displaystyle{\sum_{J \in \Gamma_{\d}^{\ell+1}} } \sigma_J^{\ell+1} B_J^{\ell+1} \,.$$
The truncation operator ${\trunc}^{\ell+1} : V^{\ell} \rightarrow V^{\ell+1}$   with respect to level $\ell+1$ acts on $s$ as follows,
$$
{\trunc}^{\ell+1} (s)   \,:=\, 
\displaystyle{ \sum_{J \in\Gamma_{\d}^{\ell+1} \,:\, 
\supp B_J^{\ell+1}\, \not\subseteq\, \Omega^{\ell+1}  } } \sigma_J^{\ell+1} B_J^{\ell+1} \,,
\quad \ell=0,\ldots,M-2.
$$ 
The cumulative truncation operator ${\Trunc}^{\ell+1} : V^{\ell} \rightarrow S_{\cal H} \subseteq V^{{M-1}}$ with respect to all finer levels in the hierarchy is then defined as 
\begin{equation*}
{\Trunc}^{\ell+1}(s) \,:=\,   {\trunc}^{M-1}( {\trunc}^{M-2}(\cdots( {\trunc}^{\ell+1}(s)) \cdots))\,,
\quad \quad \ell=0,\ldots,M-2.
\end{equation*}
Truncated hierarchical B-splines (THB-splines) are defined by exploiting the truncation mechanism \cite{giannelli2012}.
\begin{dfn}\label{dfn:thb}
The truncated hierarchical B-spline basis ${\cal T}_{\d} ({\cal G}_{\cal H})$ of degree $\d$  with respect to the mesh ${\cal G}_{\cal H}$ is defined as
\begin{equation*}  \label{truncbasis}
{\cal T}_{\d} ({\cal G}_{\cal H}) \,:=\,
%\bigcup_{l=0}^{M-1} 
\left\{
T_J^{\ell} \,:\, J \in A_\d^{\ell}, \,\, \ell=0,...,M-1  \right\}\,, \quad \text{with}\quad T_J^{\ell} \,:=\, {\Trunc}^{\ell+1} (B_J^\ell)\,,  
\end{equation*}
and $ B_J^\ell$ is called the \emph{mother} B-spline of $T_J^\ell$.
 \end{dfn}
 
THB-splines define a strongly stable basis for the hierarchical spline space $S_{\cal H}$ with respect to the supremum norm \cite{giannelli2014}. This means that there exist two constants not dependent on the number of hierarchical levels that can be multiplied to the norm of the coefficients associated to the THB-spline representation of a function $f$ to bound from below and above the norm of $f$ itself. This is not the case for the classical hierarchical basis of Definition~\ref{dfn:hb} which is only weakly stable (the associated stability constants have at most a polynomial growth
in the number of hierarchical levels) \cite{kraft1998}.

In addition, THB-splines form a convex partition of unity and preserve the coefficients of their mother functions (they preserve the coefficients of functions represented with respect to one of the bases ${\cal B}_\d^{\ell}$) \cite{giannelli2014}. This last property is indicated as \emph{preservation of coefficients} and enables the construction of {\it quasi-interpolation} operators in hierarchical spline spaces that do not require additional computations \cite{speleers2016}. Let us clarify this important property.
We denote by $V$ a tensor--product spline space of degree $\d$ defined on $\Omega^0$ and generated by a B-spline basis ${\cal B}_\d := \{B_J\,, J \in \Gamma_\d\} $ defined over a certain grid ${\cal G}$. Let $Q : S(\Omega^0) \rightarrow V$ be a quasi--interpolation operator of the form
\begin{equation}\label{scatteredqi1}
\quad Q(f) :=\sum_{J \in \Gamma_{\d}} \lambda_J(f) \, B_J\,, 
\end{equation}
for a suitable space $S$. Each functional $\lambda_J : S(\Omega^0) \rightarrow \RR$ is locally defined  and can be of discrete or continuous type. Thus, if $Q^\ell(f)$ denotes the instance of $Q$ when $V= V^\ell$ and $\lambda_J^\ell$ the related functionals for $\ell=0,\ldots,M-1$, in virtue of the preservation of coefficients, the hierarchical quasi-interpolant $Q_{\cal H} : S(\Omega^0) \rightarrow S_{\cal H} $ is simply defined as 
\begin{equation}\label{scatteredqiH} 
\quad  Q_{\cal H}(f) := \sum_{\ell=0}^{M-1} \sum_{J\in A^\ell_\d} \lambda_J^\ell(f)\, T_J^\ell.
\end{equation}
Note that $Q_{\cal H}$ has the capability of reproducing polynomials of degree $\d$, whenever $Q$ has this property. Under mild restrictions on the domain hierarchy, it is also possible to construct operators $Q_{\cal H}$ which are projectors onto the hierarchical spline space \cite{speleers2016}.
 Adaptive approximations  by different discrete THB-splines quasi-interpolation schemes  were recently introduced in the literature \cite{bracco2016, speleers2017,speleers2016}, by starting from information on $f$ (and its derivatives in the case of Hermite spline operators) available on a set of gridded points. 
 
The suitable treatment of a {\it scattered} data set, 
\begin{equation} \label{scatset}
F=\{({\bf X}_i,f_i),\, i=1,...,n,\,{\bf X}_i \in \Omega, \mbox{with } {\bf X}_i \ne {\bf X}_j \mbox{ if } i \ne j\},
\end{equation}
requires the design of adaptive schemes able to locally tailor the nature of the fitting method according to the available number of local scattered  data points. In order to deal with data of non-trivial complexity, flexible adaptive tools that generate compact representations need to be developed. This requires a versatile definition of local (data dependent) approximations and a low cost construction of the global approximation. We address these issues by: 1) computing local polynomial approximations for defining the operator $Q$ in (\ref{scatteredqi1});  2) extending them to adaptive spline spaces via hierarchical quasi-interpolation. The next two sections detail points 1) and 2).

%%%%%%%%%%%%%%%%%%%%%%%%%%%
\section{Data-dependent local polynomial least squares approximation} \label{sec:scatteredqi} 

In this section we introduce the algorithm  for the definition of the quasi-interpolant $Q$  in (\ref{scatteredqi1}) by computing the value of any functional $\lambda_J(f), J \in \Gamma_{\d},$ associated to the B-spline $B_J$ for a certain scattered data set $F$ of the form (\ref{scatset}). 

We look for a suitable local approximation in the space $\Pi_J$ of local $r$-variate polynomial of total degree $d_J \le d,$ with 
$$d := \min_{h=1,\ldots,r} d_h.$$ 
For the representation of this local polynomial, we consider the standard (local) power basis suitably translated and scaled in order to map the support of $B_J$ into the hyper-rectangle $[0,1]^r$.
Obviously, this approximation needs to be successively converted in the local spline space  $V_J$ of degree $\d$ associated with the grid ${\cal G} \cap \mbox{ supp} (B_J)$ in order to set $\lambda_J(f)$ equal to its $J$--th coefficient in the local B--spline basis. 
By following  \cite{davydov2004}, the polynomial  is obtained by  least squares approximation of a suitable local subset of data 
$$F_J \subset F,$$
initially defined as the set of data $ ({\bf X}_i,f_i)\in F \mbox{ with } {\bf X}_i $ belonging to the ball of radius $\rho_J$ equal to the halved diameter of the support of $B_J$ and with the same center. 
If $F_J$ is empty, we enlarge this local subset of data by iteratively replacing $\rho_J$ with $k \rho_J,\, k\in \NN,$ with $1< k \le K_J,$ where $K_J$ denotes a positive integer  upper bound given in input.\footnote{The upper bound $K_J$ for $k$ could cause a failure of the algorithm for some $J$  if the data have a  very highly varying density and the mesh ${\cal G}$  is not suitably chosen. In this case, a better choice of ${\cal G}$ is suggested. However, this upper bound is reasonable, since we look for a local approximation.} Once $F_J$ is fixed, we initialize $d_J = d$ and check if the minimal singular value of the collocation matrix associated with the considered local least squares problem is greater or equal than a prescribed threshold $\sigma,$ with $0 < \sigma \le 1.$  If this is not the case, $d_j$ is decreased by one. Note that, at least when $d_J$ becomes zero, the minimal singular value is surely greater or equal to $\sigma,$ since $\sigma \le 1$ and the coefficient matrix reduces to a column vector with all unit entries. The use of the threshold $\sigma$ not only guarantees that the local least squares approximation has a unique solution, but it also ensures that the quality of the polynomial  of degree $d_J$ is comparable with the one of the local best approximation in the same degree polynomial space \cite{davydov2002}.

For facilitating the comprehension of Algorithm 1, which details the procedure described above, we report %detail 
in the two following lists the considered notation.
%\newpage

\medskip
\noindent {\bf Input parameters}
\begin{itemize}
\item $F$:  full set of $n$ scattered data points, defined as in (\ref{scatset});
\item ${\cal G}$: tensor-product mesh associated with the considered global spline space $V$;
\item $\sigma \, (\le 1) $: positive threshold prescribing a lower bound for the minimal singular value (MSV) of the local least squares matrix; 
\item $K_J \, (\ge 1) $: integer upper bound which controls the maximum size of the local data sample used to construct the local least squares polynomial approximation. 
\end{itemize}

\newpage
\noindent 
% Variables and 
{\bf Further notation used in Algorithm 1}
\begin{itemize}
\item $F_J$: local set of scattered data points;
\item $I_J$: local set of indices  between $1$ and $n$ such that $F_J = \{(X_i,f_i) \in F : i \in I_J\};$ 
\item $\vert F_J \vert$: cardinality of $F_J;$
\item $d_J$:  total polynomial degree; 
\item $\Pi_J$: local $r$--variate polynomial space of total degree $d_J;$
 \item $\vert \Pi_J \vert$: dimension of $\Pi_J;$
\item ${\cal P}_J$: local power basis of $\Pi_J;$
\item $A_J$: local least squares matrix collocating ${\cal P}_J$ at the vertices of $F_J;$
\item $\mbox{MSV}(A_J)$: minimal singular value of  $A_J$;
\item $V_J$: local spline space spanned by the B-splines $B_I, I \in \Gamma_\d $, whose support intersects $\mbox{ supp} (B_J)$.
 \end{itemize}

Finally, the short notation $\lambda_J$ replaces $\lambda_J(f)$  to denote the $J$--th functional in the remaining part of the paper.
\bigskip

\noindent {\bf Algorithm 1} 
\begin{enumerate}
\item[] {\bf Input:}
\begin{itemize}
\item $F\,, {\cal G}\,, \d = (d_1,\ldots,d_r)\,, J\in \Gamma_\d$, $ \sigma \le 1$, $K_J \ge 1;$ 
\end{itemize}
\item set $\displaystyle{d := \min_{h=1,\ldots,r} d_h} ;$
\item initialize $d_J = d$ and correspondingly $\Pi_J$ and ${\cal P}_J;$
 \item set $\rho_J := (\diam(\supp B_J))/2$;
\item set $C_J $ equal to the center of $  \supp B_J $;
 \item initialize $I_J = \Big\{1\le i \le n \,:\, \,\Vert {\bf X}_i - C_J \Vert \le \,\rho_J \Big\}$ and consequently $F_J$;
 \item initialize $k=2$;
\item while  $\vert F_J \vert = 0$ and $k \le K_J$  
\begin{enumerate}
       \item  update $I_J = \Big\{i:\,\Vert {\bf X}_i - C_J \Vert \le k\,\rho_J \Big\}$  and consequently $F_J$;
       \item increase $k$ by one;
 \end{enumerate}
 \item if $\vert F_J \vert = 0$ and $k > K_J$ FAILURE 
 \item compute the matrix $A_J = [p(X_i)]_{p \in {\cal P}_J, i\in I_J}$;
\item if $\vert F_J \vert<\vert {\Pi}_J \vert$ or $\mbox{MSV}(A_J) < \sigma$ then 
\begin{enumerate}
\item decrease $d_J$ by $1$ and update $\Pi_J$ and ${\cal P}_J;$
\item update the matrix $A_J = [p(X_i)]_{p \in {\cal P}_J, i \in I_J}$;
\end{enumerate}
\item compute the least squares approximation $p_J \in \Pi_J$ of the data in $F_J$; 
\item represent $p_J$ as a function of the local spline space $V_J, \,\, p_J  = \displaystyle{ \sum_{ 
B_I \in V_J} \alpha_I B_I};$
\item set $\lambda_J = \alpha_J.$  
\item[] {\bf Output:} 
\begin{itemize}
\item the coefficient $\lambda_J$.
\end{itemize}
\end{enumerate}
\bigskip

We should mention that the possibility of using directly $V_J$  for determining a local spline approximation instead of a polynomial one has been discarded since not feasible in general. The motivation is that $V_J$ is a suitable choice only if the projection into $\Omega$ %$\Omega^0$ 
of $F_J$ is not larger than the support of $B_J$. In this case,  step 7  of the algorithm should be removed and a FAILURE of the algorithm is more likely. As confirmed by the numerical experiments, for data sets with  highly varying density, this can  compromise the benefit of the hierarchical formulation in (\ref{scatteredqiH}) of $Q,$ since very few levels and relatively coarse meshes can be used.  On the other hand, the choice of a polynomial space $\Pi_J$ %$P_J$ 
with an adaptively selected total degree $d_J \le d,$ together with step 7  when a suitable $K_J$ is considered, guarantees the success of the hierarchical implementation described in the next section (under the mild requirement of a suitable choice of ${\cal G}^0$ --- see also Remarks \ref{first} and \ref{inter} below). Note that the numerical tests have also confirmed that for some real world data sets, replacing in the algorithm tensor-product polynomial spaces to total degree ones can produce artifacts. Since it would also be more expensive, we consider only polynomial spaces of total degree.

%---------------------------------------------------------------------------
\section{Adaptive scattered data quasi-interpolation} \label{sec:adap}

\medskip
If the subdomain hierarchy is given, the quasi--interpolation operator $Q$ introduced in the previous section and based on Algorithm 1 can be easily extended to a hierarchical spline space using the definition in (\ref{scatteredqiH}). However, in order to get an effective and fully automatic approximation scheme able to control the maximum approximation error at the data points, while reducing the number of degrees of freedom, a strategy for adaptively defining $M$ and ${\cal G_H}$ is necessary. This is the goal of Algorithm 2 below which successively modifies the number of levels and the hierarchical mesh by comparing the current values of the errors with a given suitable positive tolerance $\epsilon.$ Whenever ${\cal G_H}$ is updated, the quasi--interpolant $Q_{\cal H}(f)$ needs to be updated. Algorithm 1 is then called for  the computation of each new functional $\lambda_J^\ell$ with $\ell$ varying between $1$ and the current number of levels. If ${\cal G}^0$ is suitably chosen and the maximum number of levels $M_{max}$ is sufficiently high, the algorithm fails only if a too small tolerance $\epsilon$ is given in input.
 \noindent
Before presenting the algorithm,  we need to introduce some further notation,
\begin{itemize}
\item $\rho_J^\ell := \diam(\supp B_J^\ell)/2$, for $J \in \Gamma_\d^\ell, \,\ell=0,...,M-1;$
\item $\delta^\ell := \max_{J\in \Gamma_\d^\ell}\rho_J^\ell$, for $\ell=0,...,M-1.$
\end{itemize}
In order to simplify the notation, we extend the two previous definitions also to $\ell=-1$, by simply considering a B-spline set ${\cal B}_\d^{-1}:=\{ B_J^{-1}, \,J\in \Gamma_\d^{-1}\}$ defined on an auxiliary tensor-product mesh whose size is two times the size of ${\cal G}^0$ and having at least $ (d_1+1)\cdots(d_r+1)$ 
elements.

\begin{figure}[!th]%\begin{center}%.31
\centering{
\hspace*{-0.6cm}
\subfigure[]{{\includegraphics[trim=20mm 10mm 25mm 10mm,clip,scale=0.31]{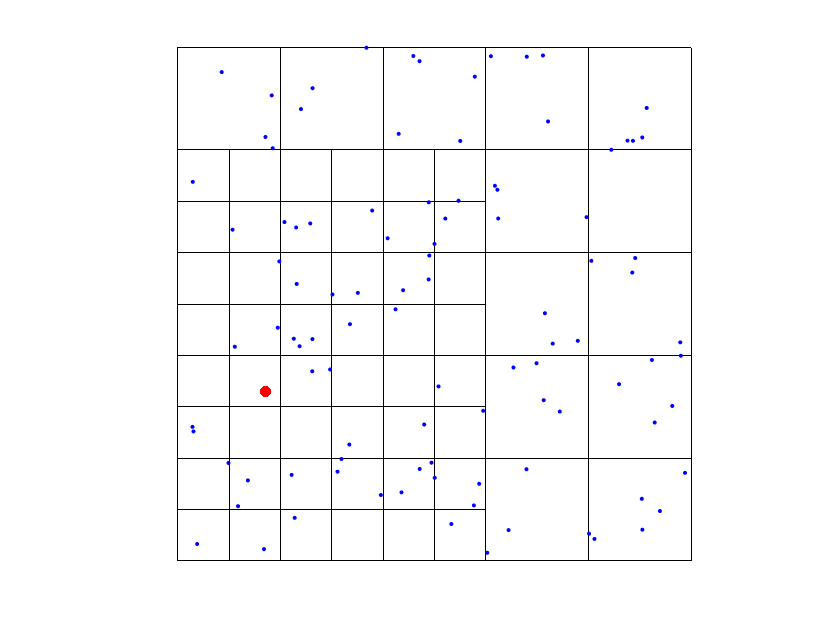}}}\hspace*{-0.3cm}
\subfigure[]{{\includegraphics[trim=25mm 10mm 25mm 10mm,clip,scale=0.31]{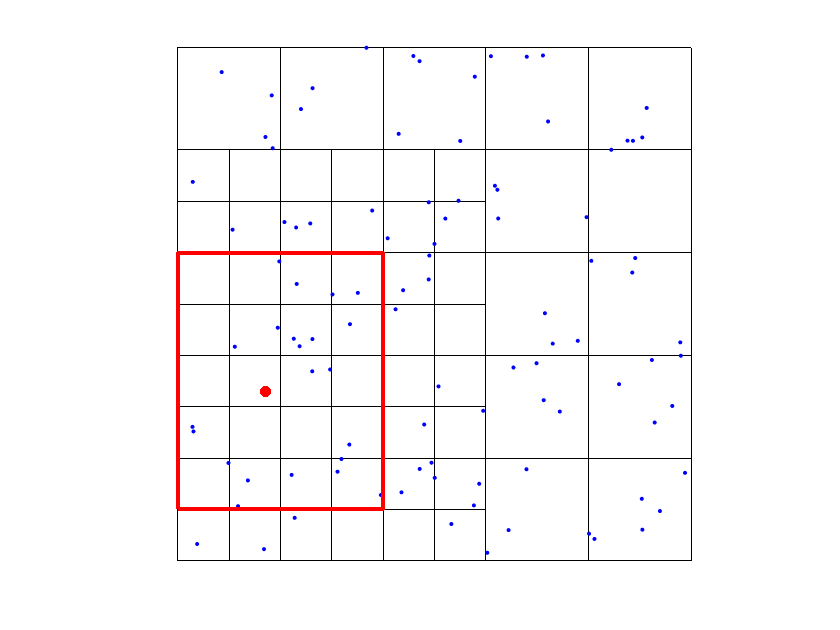}}}\hspace*{-0.3cm}
\subfigure[]{{\includegraphics[trim=25mm 10mm 15mm 10mm,clip,scale=0.31]{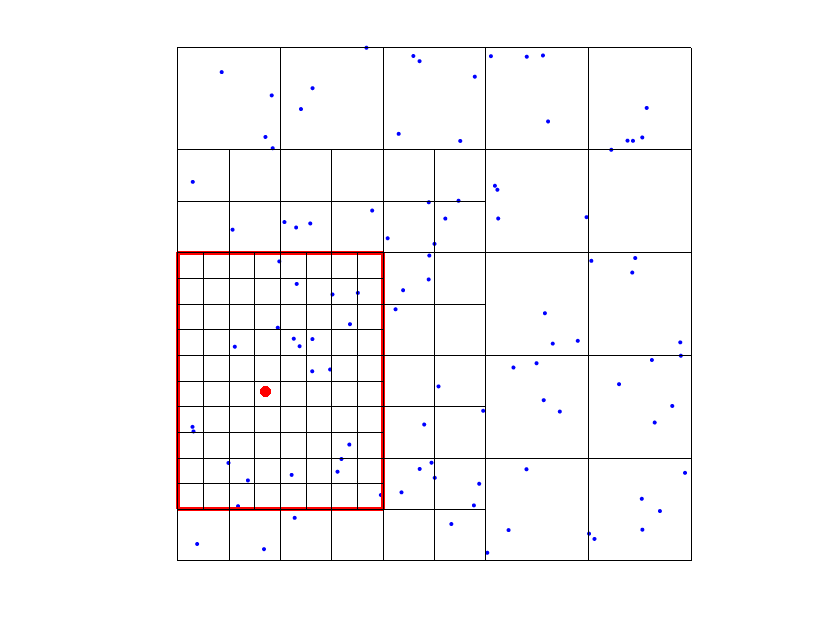}}}
}
\caption{Example of refinement with $\d=(2,2)$: for each data point that does not satisfy the given  tolerance (e.g. the red point in (a)), we mark the union of the B-spline supports which contain it (highlighted in (b)), and then subdivide the related cells (c).}
 \label{funref}
%\end{center}
\end{figure}

In order to define the refined hierarchical mesh in the algorithm, we use a refinement by functions: for each data point where the tolerance is not satisfied, we select all the B-splines $B_J^\ell$ such that $T_J^\ell\in {\cal T}_{\d} ({\cal G}_{\cal H})$ and whose support contains the data point. Then, we subdivide all the active cells (that belong to ${\cal G}_{\cal H}$) of the supports of the selected functions (see Figure \ref{funref}). Note that refining the supports, by Definitions \ref{dfn:hb} and \ref{dfn:thb}, implies that we are always replacing some functions in the THB-splines basis ${\cal T}_\d({\cal G}_{\cal H})$ with more functions on a finer level. In other words, this approach guarantees to add degrees of freedom in the areas where we need a further error reduction.

\bigskip
\noindent {\bf Algorithm 2} 
\begin{enumerate}
\item[] {\bf Input:}
\begin{itemize}
\item $F\,, {\cal G}^0\,, \d = (d_1,\ldots,d_r)\,; $ 
\item maximum number of levels $M_{max} \ge 2$ and tolerance $\epsilon > 0$ for adaptive refinement;
\item positive threshold $ \sigma \le 1$ governing  the degree selection for each local polynomial approximation.
\end{itemize}
\item set ${\cal G}_{\cal H}={\cal G}^0, S_{\cal H} = V^0$ and $M=1$;
\item for each $J \in \Gamma_{\d}^0$,  compute $\lambda_J^0$ with Algorithm 1 using as input $F, \,{\cal G} = {\cal G}^0, \, \d\,, J \in \Gamma_{\d}^0, \,\sigma$ and
\begin{equation*}
 K_J :=\,  \Bigg(\Bigg\lceil
   \frac{2\,\delta^{-1}}{\rho_J^0}
 \Bigg \rceil +1\Bigg);
\end{equation*}
\item if for some $\lambda_J^0$ Algorithm 1 ends with a FAILURE, end with a FAILURE;
\item   compute the errors
\begin{equation*}
e_i:=\vert Q_{\cal H}(f)({\bf X}_i)-f_i\vert, \qquad i=1,...,n;
\end{equation*}
\item  {while} $\max_{i=1,...,n}e_i>\epsilon$ and $M \le M_{max}$

\begin{enumerate}
\item  {for} $i=1,...,n$, if $e_i>\epsilon$,  mark the functions $B_J^{\ell},\,  J\in A_{\d}^\ell\,, \, \ell=0,\ldots,M-1$   such that ${\bf X}_i\in\supp B_J^{\ell}$;
\item  define the new mesh ${\cal G_H}$ and the corresponding space $S_{\cal H}$ by performing dyadic refinement of the cells belonging to the support of the marked functions;
\item if ${\cal G_H}$ includes cells belonging to ${\cal G}^{M+1},$ increase $M$ by one;
\item if $M > M_{max},$ FAILURE; 
\item  {for} $\ell=0,...,M-1$ and for each new $J \in A_{\d}^\ell$, compute $\lambda_J^{\ell}$ with Algorithm 1 using as input $F, \,{\cal G} = {\cal G}^\ell, \, \d\,, J \in A_{\d}^\ell, \,\sigma$ and
\begin{equation*}
 K_J := \Bigg(\Bigg\lceil\frac{2\delta^{\ell-1}}{\rho_J^\ell}\Bigg\rceil+1\Bigg);
\end{equation*}
\item update the errors
\begin{equation*}
e_i:=\vert Q_{\cal H}(f)({\bf X}_i)-f_i\vert, \qquad i=1,...,n.
\end{equation*}
\end{enumerate}

\item[] {\bf Output:} 
\begin{itemize}
\item $M \le M_{max}$ and the hierarchical mesh ${\cal G}_{\cal H}$;
\item all the coefficients $\lambda_J^\ell, J \in A_\d^\ell, \ell =0,\ldots,M-1 $ of the hierarchical quasi-interpolant $Q_{\cal H}(f)$. 
\end{itemize}
\end{enumerate}

\noindent
The second possible FAILURE of the algorithm in step 5(d) could sometimes be removed by simply increasing $M_{max}.$ However it could also persist if $\epsilon$ is too small (but clearly this  tolerance can not be arbitrarily small).  The following two remarks discuss two important aspects of the algorithm. Remark \ref{first} focuses on the first possible FAILURE of the algorithm and on how it is possible to avoid it. Remark \ref{inter} explains why, if step 2 is overcome, then in the ``{while}'' cycle the required new functionals can always be computed.

\begin{remark} \label{first}
Because of the initialization of $K_J$ in step 2, we cannot guarantee, in general, that it is possible to compute all the coefficients $\lambda_J^0$ in this step. This motivates step 3 of the algorithm. This initial FAILURE could be avoided by initializing $K_J$ with a sufficiently high integer value. However, this is not the right solution, since this would require to consider data $({\bf X}_i,f_i)$ with ${\bf X}_i$ quite far from the supports of the  THB-splines. Our experiments confirm that this can significantly worsen (and sometimes also ruin) the accuracy of the quasi-interpolant. In order to avoid this problem, we suggest a more careful choice of ${\cal G}^0,$ which takes into account the chosen degree as well as the data distribution. 
\end{remark}

\begin{figure}[!h]\begin{center}
\subfigure[]{{\includegraphics[trim=20mm 10mm 20mm 10mm,clip,scale=0.25]{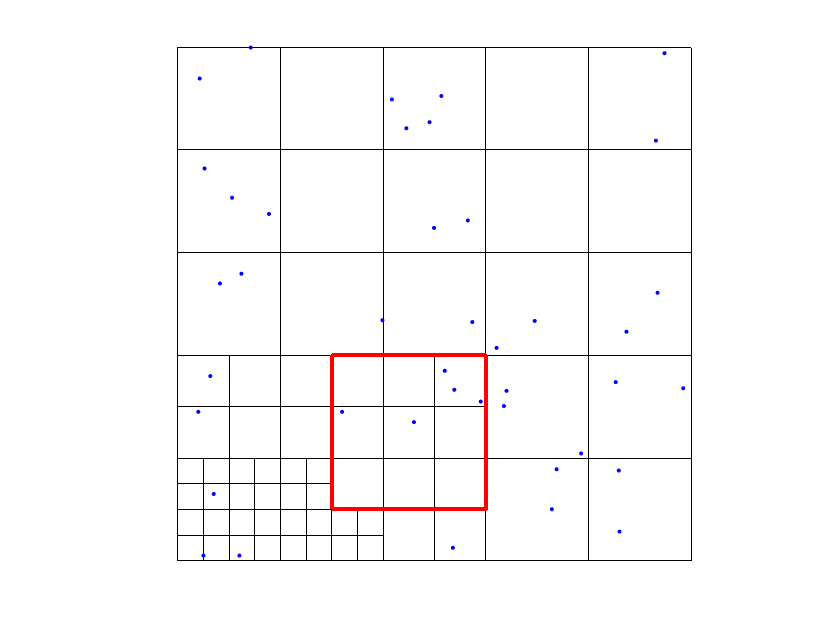}}}%\hspace*{0.3cm}
\subfigure[]{{\includegraphics[trim=20mm 10mm 20mm 10mm,clip,scale=0.25]{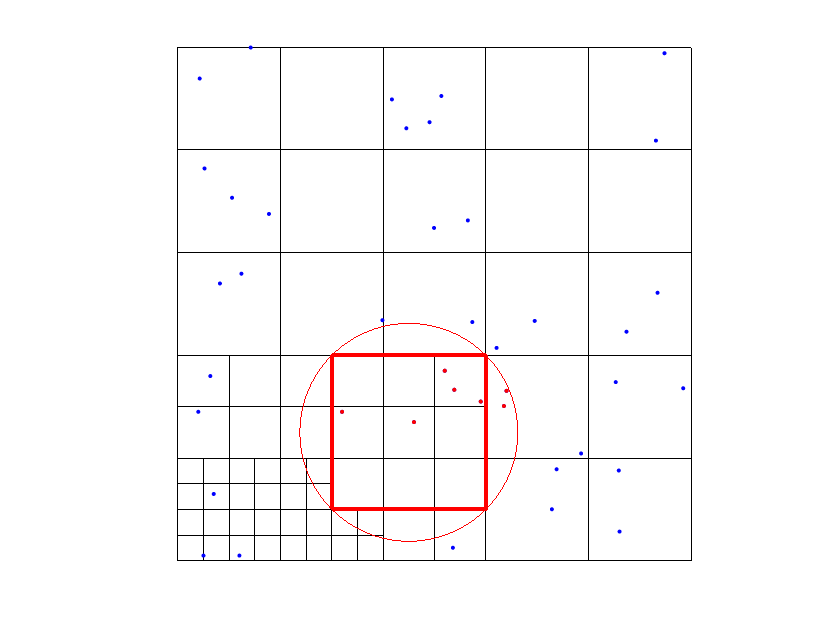}}}\\
\subfigure[]{{\includegraphics[trim=20mm 10mm 20mm 10mm,clip,scale=0.25]{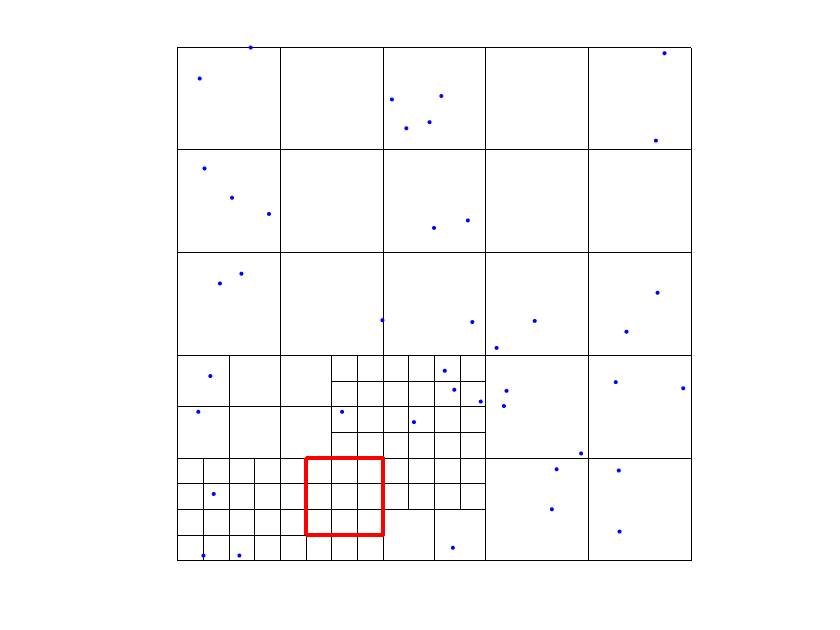}}}
\subfigure[]{{\includegraphics[trim=20mm 10mm 20mm 10mm,clip,scale=0.25]{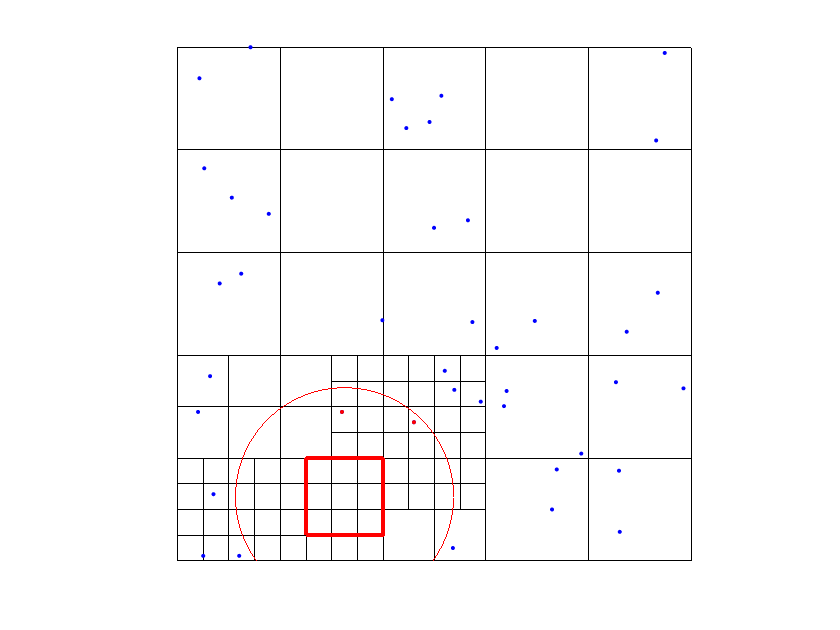}}}
\caption{For $\d=(2,2),$ the support of a B-spline $B_J^\ell$ to be refined is highlighted in (a) while the circle and the points used to determine $\lambda_J^\ell$ are shown (red) in (b). The support of a new B--spline $B_I^{\ell+1}$ intersecting the support of $B_J^\ell$ is highlighted in (c) and, according to step 5(e) of Algorithm 2 and to step 7 of Algorithm 1, the related circle and points used to determine $\lambda_I^{\ell+1}$ are shown in (d).}
 \label{locref}
\end{center}
\end{figure}

\begin{remark} \label{inter}
The choice of $K_J$ addresses several issues. First, it was designed to keep it strictly related to the size of the mesh of the level $\ell$ of the considered $T_J^\ell.$ In fact, as mentioned in the previous remark, our tests suggests that considering data points too far from the supports of the THB-splines often worsen the accuracy of the quasi-interpolant. Moreover, our setting in step 5(e) guarantees that, at each ``{while}'' cycle  we have locally enough data points to compute the coefficients associated with the new active $T_J^\ell.$ First, consider that the new active THB-splines are of levels $\ell>0.$ Note then that for the computation of the new $\lambda_J^\ell$, $\ell>0$, we look for data points at most in a ball which includes the supports of the marked functions intersecting the support of $B_J^\ell$. Since by Step 4(a) the support of the previously marked functions contains at least one data point, this guarantees that we have locally at least one data point to compute $\lambda_J^\ell$ (see Figure \ref{locref}).  
\end{remark}

%%%%%%%%%%%%%%%%%%%%%%%%%%%%%%%%%%%%%%%%%%%%%%%%%%%%%%%%%%%%%%%%%%%%
  %%%%%%%%%%%%%%%%%%%%%%%%%%%%%%%%%%%%%%%%%%%%%%%%%%%%%%%%%%%%%%%%%%%%
\section{Numerical results}\label{sec:exm}

In this section we present the results obtained by applying the hierarchical approximation algorithm to different types of data, including comparisons with other spline approximation methods, when available.

For the tests presented in Examples 1-4, we always use $\d=(2,2)$.
If not stated otherwise, we always choose $M_{max}$ sufficiently large to be able to meet  the given tolerance $\epsilon$, which is specified in the corresponding tables (Tables \ref{test1a},\ref{test2a},\ref{test3a},\ref{test4a}).

In the following tests, we report the number of degrees of freedom (NDOF) and the two errors
\begin{equation*}
e_{max}:=\max_{1\le i\le n} e_i, \qquad e_{RMS}:=\sqrt{ \frac{1}{n} \displaystyle{\sum_{i=1}^n} e_i^2}.
\end{equation*}
In all the cases, the results show that the algorithm succeeds at providing accurate approximations of the data satisfying the given tolerance on the maximum error. In Examples \ref{ex1}-\ref{ex2}, we also showed that the method favourably  compares with respect to previous two-stage methods that consider spline spaces over uniform triangulations. In particular, with our adaptive scheme we get essentially the same maximum error with significantly less degrees of freedom. We note that the $e_{RMS}$ is instead slightly higher, which is probably natural since the refinement in Algorithm 2 is completely driven by a control based only on the maximum error $e_{max}$.

\begin{exm}\label{ex1}
In the first example we approximate the black forest elevation data set ($15885$ points) already used for testing scattered data approximation methods in \cite{davydov2004} and \cite{davydov2006}, where two-stage methods with splines over uniform triangulations were considered. This is a case of scattered data with highly varying density, as evident from Figure \ref{bf1}(a). Figure \ref{bf1}(b) clearly shows the adaptive nature of our method, which produces a very refined mesh in the areas of the domain where the data density is high and it is harder to get the required accuracy (they correspond to mountain regions, as it can be seen in Figures \ref{bf2}(a)-(c)). In Table \ref{test1a}, we report the step-by-step results obtained with our method.  

In Table \ref{test1b} we compare the performance with two approximants on scattered data defined in spline spaces over triangulations: the first, denoted by $P$, is obtained in \cite{davydov2004} by extending local polynomial approximation, while the second, denoted by $H_{MQ}$ is constructed in \cite{davydov2006} by using local hybrid polynomial/radial basis approximations. %In Table \ref{test1c} we show how the number of data employed to compute the local approximations (and then the coefficients of $Q_{\cal H}$) varies depending on the level $\ell$.

\begin{table}[H]
\centerline{ 
\begin{footnotesize}\begin{tabular}{SSSSS} \toprule
		{$M$} & {elements of ${\cal G}^{M-1}$} & {NDOF} & {$e_{max}$} & {$e_{RMS}$} \\ \midrule
1  & {$32\times 32$}         &    {1156}  & {2.680e+02}   & {6.839e+01}  \\
2  & {$64\times 64$}         &   {4058}  & {1.803e+02} & {3.630e+01}  \\ 
3  & {$128\times 128$}     & {12705}  & {7.165e+01} & {1.240e+01} \\
4  & {$256\times 256$}     & {20564}  & {5.349e+01} & {8.855e+00} \\
5  & {$512\times 512$}     & {22759}  & {3.438e+01} & {8.637e+00} \\
6  & {$1024\times 1024$} & {23092}  & {2.979e+01} & {8.620e+00} \\
 \bottomrule
\end{tabular}\end{footnotesize}
}
\caption{Numerical results for Example $1$ (black forest data set) with tolerance $\epsilon=30.0$ and $\sigma = 5\cdot 10^{-2}$.}\label{test1a}
\end{table}

\begin{table}[H]
\centerline{ 
\begin{footnotesize}\begin{tabular}{SSSS} \toprule
		{method} & {NDOF}& {$e_{max}$} & {$e_{RMS}$} \\ \midrule
	  {$Q_{\cal H}$} & {$23092$}  & {2.979e+01} & {8.620e+00}  \\
		{$P$} & {$91526$} & {3.060e+01} & {3.270e+00} \\
    {$H_{MQ}$} & {$91526$} & {3.200e+01}  & {2.170e+00} \\ \bottomrule
\end{tabular}\end{footnotesize}
}
\caption{Comparison of the performances of $Q_{\cal H}$ with the scattered data approximants $P$ and $H_{MQ}$, constructed in \cite{davydov2004} and \cite{davydov2006}, respectively, for Example \ref{ex1}.}\label{test1b}
\end{table}

\begin{figure}[!th]%\begin{center}%.4-.3
\centering{\hspace*{-1.5cm}
\subfigure[]{{\includegraphics[trim=0mm 0mm 0mm 0mm,clip,scale=0.4]{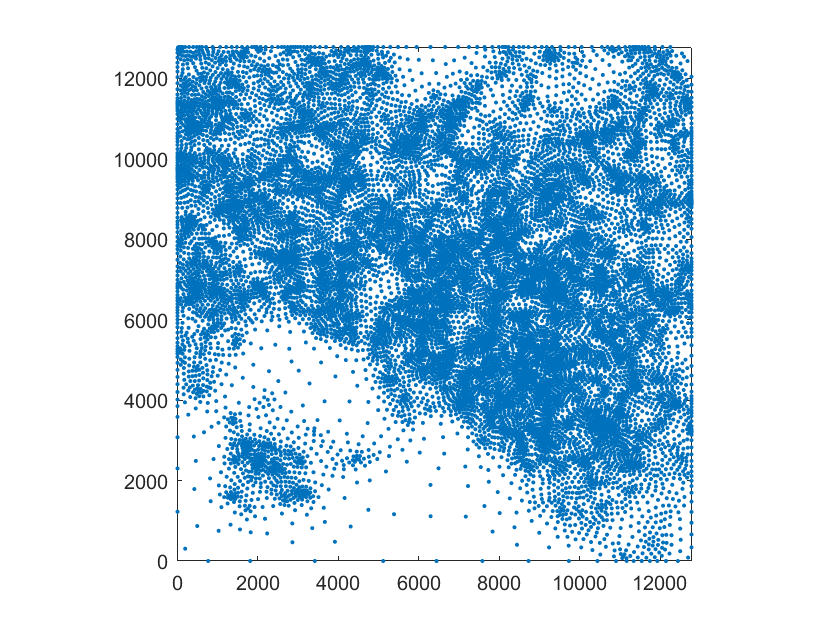}}}%\hspace*{-1cm}
\subfigure[]{\raisebox{7mm}{\includegraphics[trim=0mm 0mm 0mm 0mm,clip,scale=0.44]{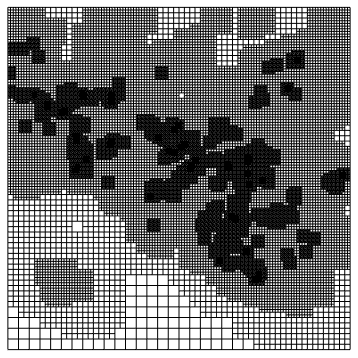}}}}
\caption{The data points of Example \ref{ex1} (a) and the hierarchical mesh generated by our approximation algorithm (b).}
 \label{bf1}
%\end{center}
\end{figure}

\begin{figure}[!th]%\begin{center}
\centering{
\subfigure[]{ 
\includegraphics[trim=0mm 30mm 0mm 35mm,clip,scale=0.55]{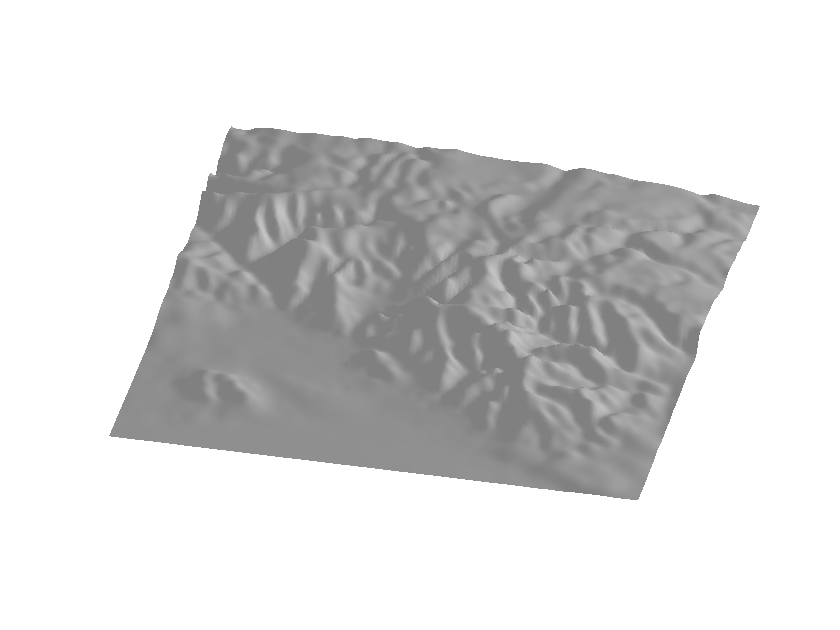}}}\\%\hspace*{0.3cm}
\subfigure[]{\includegraphics[trim=0mm 0mm 0mm 0mm,clip,scale=0.35]{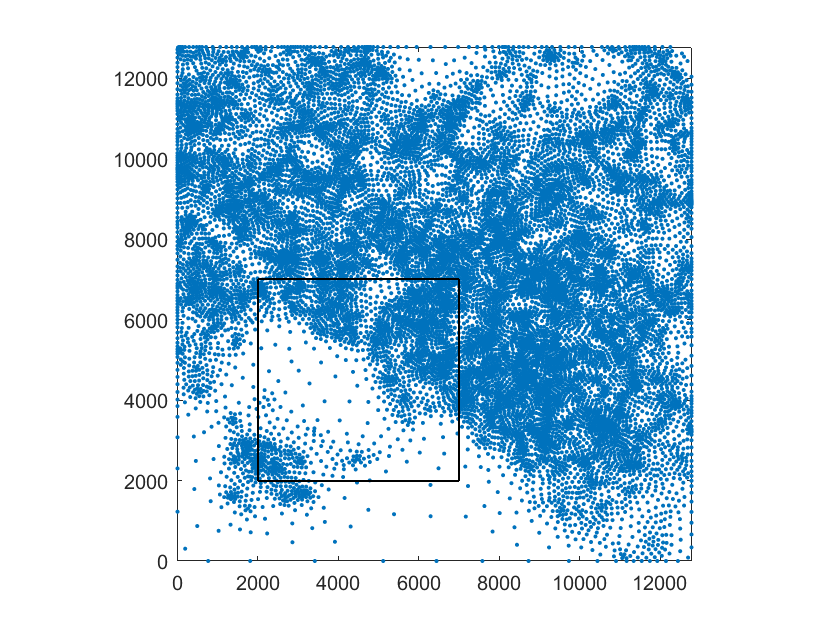}}
\hspace*{-1cm}
%\subfigure[]{\includegraphics[trim=0mm 0mm 0mm 0mm,clip,scale=0.35]{appsurfbfsection2b.png}}
\subfigure[]{\includegraphics[trim=0mm 0mm 0mm 0mm,clip,scale=0.4]{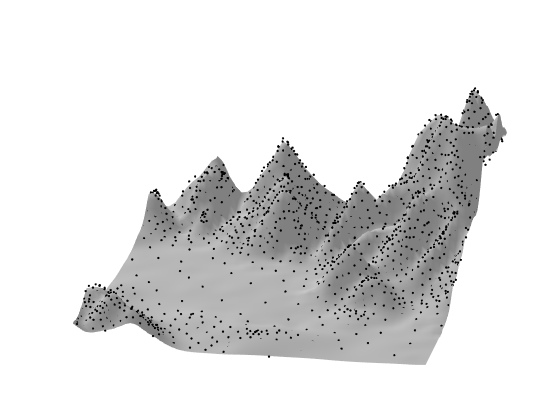}}
\caption{The resulting surface approximating the data set of Example \ref{ex1} (a); (c) is a zoom of the surface in (a), with the original data represented as black dots, corresponding to the area highlighted in (b).}
 \label{bf2}
%\end{center}
\end{figure}

\end{exm}

\begin{exm}\label{ex2}
We consider the glacier data set ($8345$ points) used for testing scattered data approximation methods in \cite{davydov2004} and \cite{davydov2006}. We here deal with a quite different type of scattered data, namely contour track data, where the data comes from the sampling along curves of equal elevation, see Figure \ref{gl2}(b). Since in this case the approximation algorithm produces a surface which, while well approximating the data, shows some unwanted oscillations, we consider a slightly different version of the scheme. In particular, we decrease the degree of the local approximations also if their evaluation at the vertices of a suitable local grid returns values too far from the range of the scattered data locally considered. Moreover, this is a case where the data points do not cover a rectangular domain: as evident from Figure \ref{gl2}(b), there are two areas (bottom left and top right corners) without any data points. We then consider a domain obtained by removing two small areas from the corners of the smallest rectangle containing all the data points, see Figure \ref{gl1}. The approximation algorithm is then applied by simply considering as $V^\ell$, $\ell=0,...,M-1$ the spline spaces defined on the rectangle containing the domain, and then discarding the B-spline functions whose support does not intersect the domain. Consequently, the hierarchical mesh in Figure \ref{gl1} shows the mesh limited to the domain. In Table \ref{test2a} and Table \ref{test2b}, we report the step-by-step results of the test and the comparison with the two approximants on scattered data defined in spline spaces over uniform triangulations, $P$ and $H_{MQ}$, provided in \cite{davydov2004} and \cite{davydov2006}, respectively. 

\begin{table}[H]
\centerline{ 
\begin{footnotesize}\begin{tabular}{SSSSS} \toprule
		{$M$} & {elements of ${\cal G}^{M-1}$} & {NDOF} & {$e_{max}$} & {$e_{RMS}$} \\ \midrule
    1  & {$16\times 16$}     & { 309}  & {5.784e+01} & {1.313e+01}  \\
    2  & {$32\times 32$}     & {1002} & {3.436e+01} & {6.381e+00}  \\
    3  & {$64\times 64$}     & {2113} & {2.217e+01} & {4.174e+00}  \\ 
    4  & {$128\times 128$} & {2626} & {1.647e+01} & {3.949e+00} \\
	5  & {$256\times 256$} & {2736} & {1.583e+01} & {3.934e+00} \\ 
 \bottomrule
\end{tabular}\end{footnotesize}
}
\caption{Numerical results for Example \ref{ex2} (glacier data set), with tolerance $\epsilon=16.0$ and $\sigma = 2\cdot 10^{-1}$.}\label{test2a}
\end{table}

\begin{table}[H]
\centerline{ 
\begin{footnotesize}\begin{tabular}{SSSS} \toprule
		{method} & {NDOF}& {$e_{max}$} & {$e_{RMS}$} \\ \midrule
	  {$Q_{\cal H}$} & {$2736$}  & {1.583e+01} & {3.934e+00}  \\
    {$P$} & {$7254$} & {1.870e+01} & {2.780e+00} \\
    {$H_{MQ}$} & {$7254$} & {1.560e+01}  & {2.260e+00}\\ \bottomrule
\end{tabular}\end{footnotesize}
}
\caption{Comparison of the performances of $Q_{\cal H}$ with the scattered data approximants $P$ and $H_{MQ}$, constructed in \cite{davydov2004} and \cite{davydov2006}, respectively, for Example \ref{ex2}.}\label{test2b}
\end{table}

\begin{figure}[!th]%\begin{center}
\centerline{\includegraphics[trim=0mm 0mm 0mm 0mm,clip,scale=0.55]{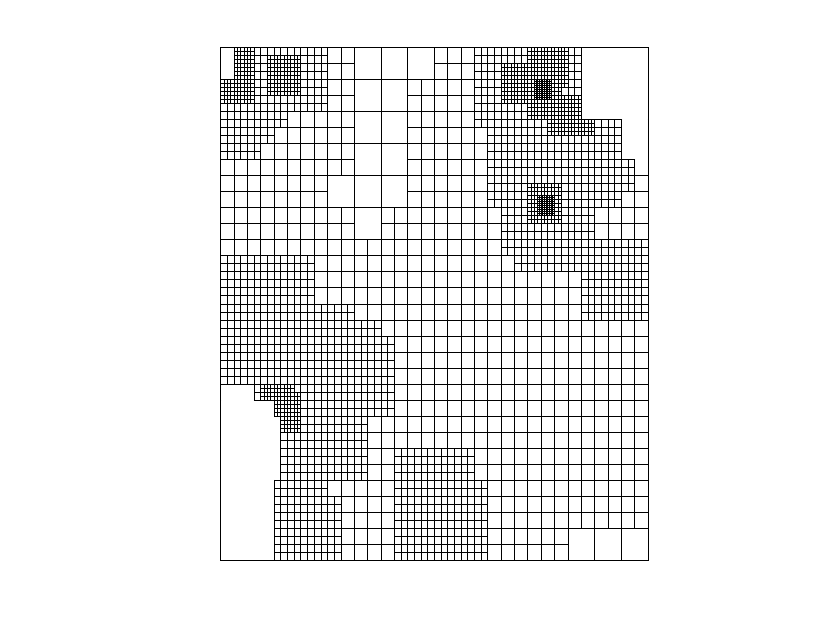}}
\caption{The hierarchical mesh generated by our approximation algorithm for Example \ref{ex2}.}
 \label{gl1}
%\end{center}
\end{figure}

\begin{figure}[!th]%\begin{center}
\centering{
\hspace*{-1.5cm}
\subfigure[]{\includegraphics[trim=0mm 40mm 0mm 35mm,clip,scale=0.85]{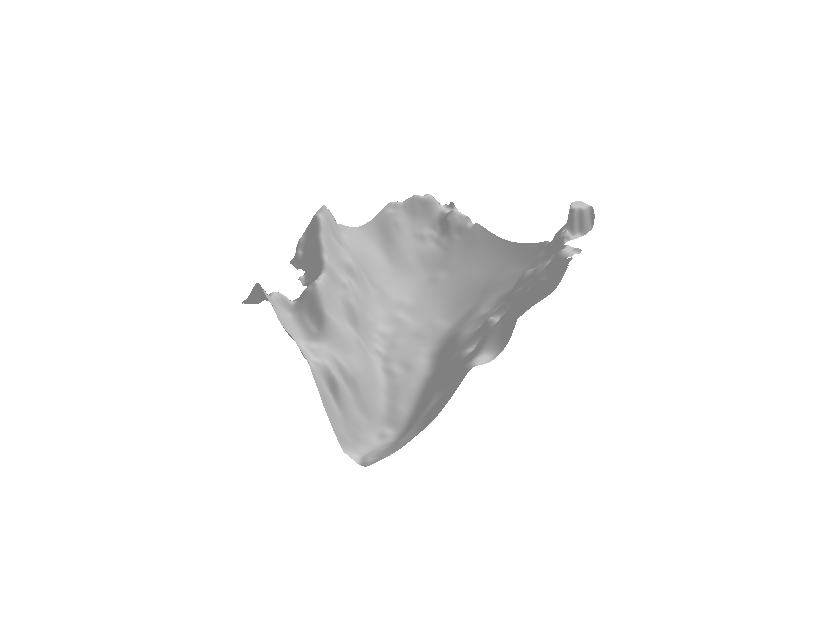}}}\\%\hspace*{0.3cm}
\subfigure[]{\includegraphics[trim=10mm 10mm 10mm 10mm,clip,scale=0.40]{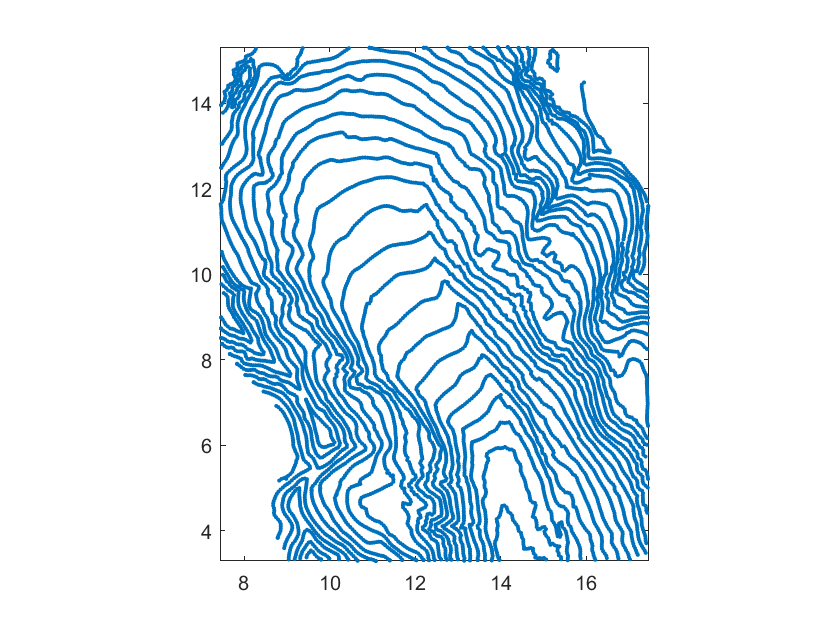}}
\hspace*{-1.5cm}
\subfigure[]{\includegraphics[trim=10mm 10mm 10mm 10mm,clip,scale=0.40]{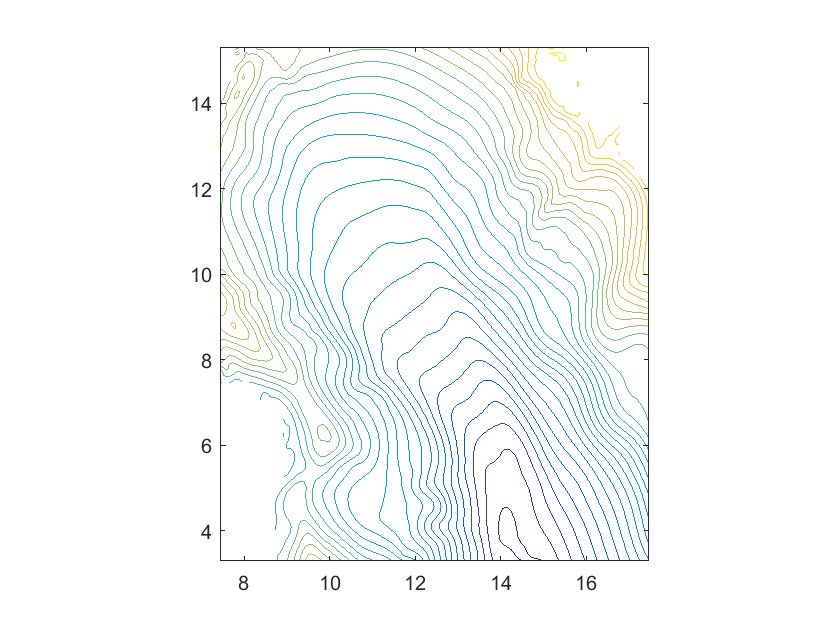}}
\caption{The resulting surface approximating the data set of Example \ref{ex2} (a); (b) and (c) are the original contour track data and the contour plot of the approximating surface, respectively.}
 \label{gl2}
%\end{center}
\end{figure}

\end{exm}

\begin{exm}\label{ex3}
We consider a subset of the Rotterdam harbour data set used in \cite{davydov2004} and \cite{davydov2006}. More precisely, we select the data whose location is inside an L-shaped domain, to give a further example of application of the algorithm to non-rectangular domains. In this case the data are affected by noise and by the presence of several outliers. In order to clean the data, we performed a preliminary step, as suggested in \cite{davydov2004}. First, we ran our approximation method on the data set using $M_{max}=6$ and tolerance $\epsilon=6\cdot 10^{-2}$. Then, we removed from the data set all the points where the error of the approximation exceeds $e_{RMS}$. The resulting data set ($12250$ points) is shown in Figure \ref{ro1}(a). Table~\ref{test3a} presents the results obtained at the different refinement steps and Figure~\ref{ro2} shows the hierarchical spline surface that approximates the given data.

\begin{table}[H]
\centerline{ 
\begin{footnotesize}\begin{tabular}{SSSSS} \toprule
		{$M$} & {elements of ${\cal G}^{M-1}$} & {NDOF} & {$e_{max}$} & {$e_{RMS}$} \\ \midrule
	1  & {$8\times 8$}             &       {91}  & {7.742e-01} & {1.631e-01}  \\
    2  & {$16\times 16$}         &     {278}  & {5.054e-01} & {1.031e-01}  \\
    3  & {$32\times 32$}         &     {941}  & {5.258e-01} & {6.677e-02}  \\
    4  & {$64\times 64$}         &   {3451}  & {3.992e-01} & {4.326e-02}  \\ 
    5  & {$128\times 128$}     & {11247}  & {3.607e-01} & {2.867e-02} \\
	6  & {$256\times 256$}     & {17018}  & {1.614e-01} & {2.498e-02} \\ 
	7  & {$512\times 512$}     & {18095}  & {9.574e-02} & {2.469e-02} \\ 
	8  & {$1024\times 1024$} & {18218}  & {6.994e-02} & {2.465e-02} \\
 \bottomrule
\end{tabular}\end{footnotesize}
}
\caption{Numerical results of Example \ref{ex3} (Rotterdam harbour data set), with tolerance $\epsilon=7\cdot 10^-2$ and $\sigma = 5\cdot 10^{-2}$.}\label{test3a}
\end{table}

\begin{figure}[!th]%\begin{center}
\hspace*{-.5cm}
\subfigure[]{{\includegraphics[trim=0mm 0mm 0mm 0mm,clip,scale=0.38]{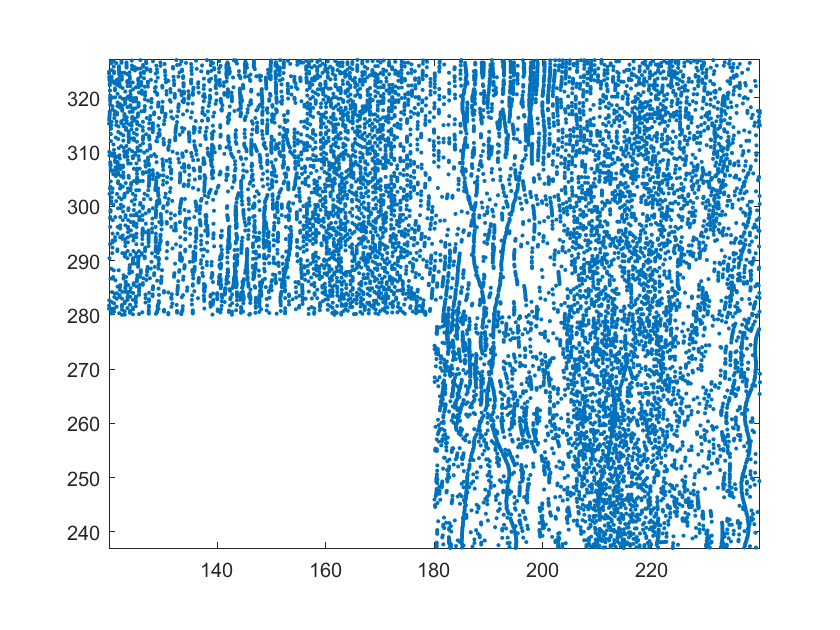}}}\hspace*{-.75cm}%\raisebox{7mm}
\subfigure[]{{\includegraphics[trim=0mm 0mm 0mm 0mm,clip,scale=0.38]{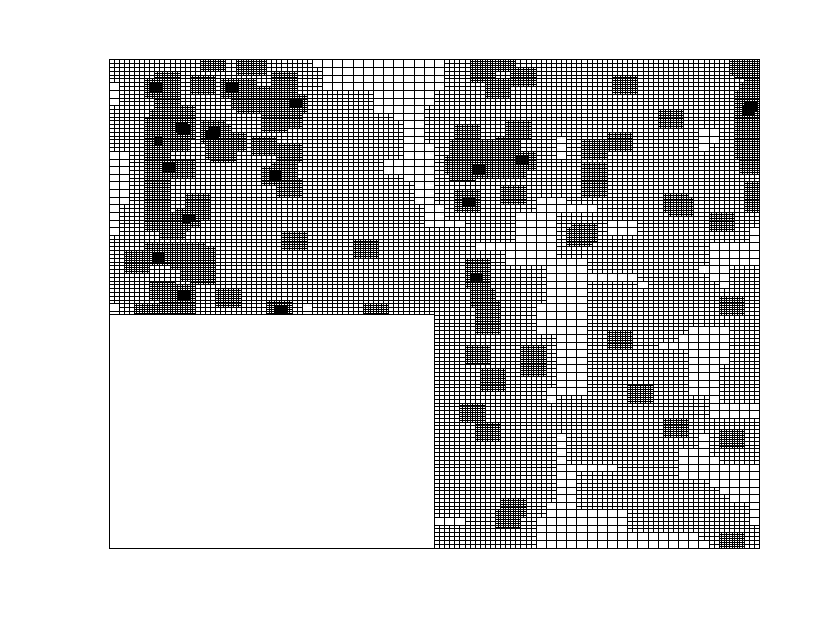}}}
\caption{The data points of Example \ref{ex3} (a) and the hierarchical mesh generated by our approximation algorithm (b).}
 \label{ro1}
%\end{center}
\end{figure}

\begin{figure}[!th]%\begin{center}
\centerline{\includegraphics[trim=30mm 40mm 30mm 25mm,clip,scale=0.70]{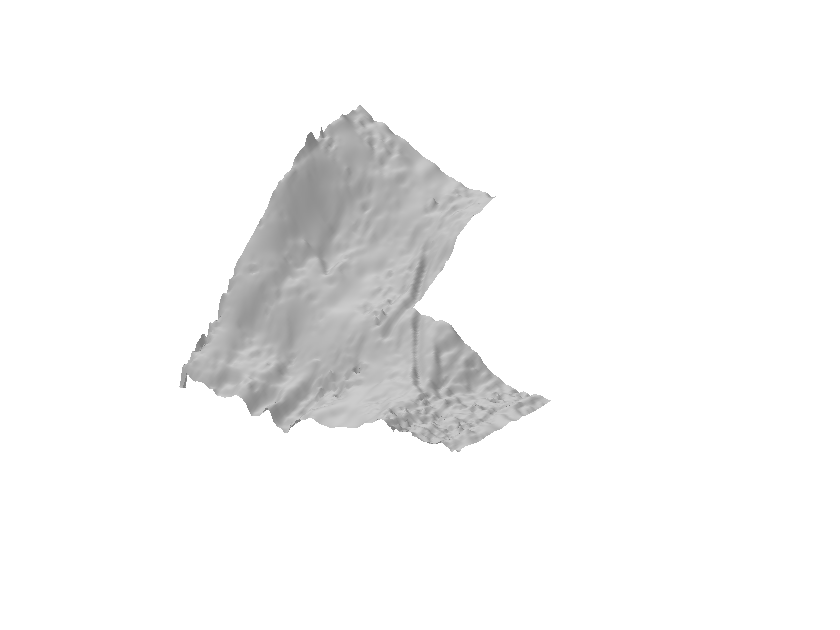}}%\hspace*{0.3cm}
\caption{The resulting surface approximating the data set of Example \ref{ex3}.}
 \label{ro2}
%\end{center}
\end{figure}

\end{exm}

\begin{exm}\label{ex4}
We present this test to assess the behaviour of our method also on gridded data. We consider a set of $129\times 129$ elevation data of a mountain region of the Hawaii Islands (available at \cite{usgeo}). In Table \ref{test4a} we report the step-by-step results of the test, while Figure~\ref{ha1} shows the hierarchical mesh obtained with the adaptive scheme. The original data surface is illustrated in Figure~\ref{ha2} together with the hierarchical spline approximation and the corresponding contour plots. Note that for exact gridded data sets other quasi-interpolation operators may be more effective.

\begin{table}[H]
\centerline{ 
\begin{footnotesize}\begin{tabular}{SSSSS} \toprule
		{$M$} & {elements of ${\cal G}^{M-1}$} & {NDOF} & {$e_{max}$} & {$e_{RMS}$} \\ \midrule
    1  & {$64\times 64$}     &   {4356} & {1.415e+02} & {1.051e+01}  \\
    2  & {$128\times 128$} & {16065}  & {7.869e+01} & {6.109e+00}  \\ 
    3  & {$256\times 256$} & {36933}  & {5.262e+01} & {4.400e+00} \\
	4  & {$512\times 512$} & {59706}  & {1.000e+01} & {2.875e+00} \\ 
 \bottomrule
\end{tabular}\end{footnotesize}
}
\caption{Numerical results for Example \ref{ex4} (Hawaii data set), with tolerance $\epsilon=10.0$ and $\sigma = 10^{-4}$.}\label{test4a}
\end{table}

\begin{figure}[h!]%\begin{center}
\centerline{{\includegraphics[trim=0mm 0mm 0mm 0mm,clip,scale=0.55]{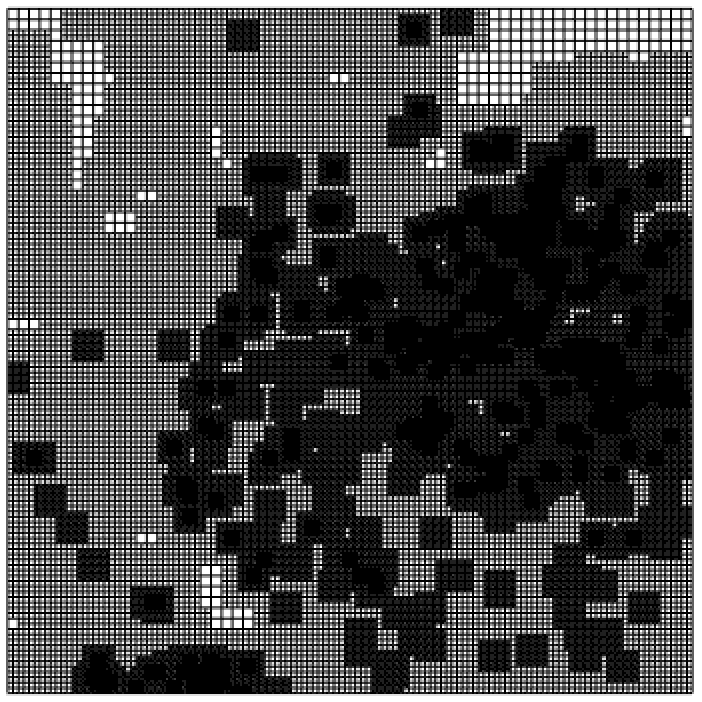}}}
\caption{Hierarchical mesh generated by our approximation algorithm in Example \ref{ex4}.}
 \label{ha1}
%\end{center}
\end{figure}

\begin{figure}[!th]%\begin{center}
\hspace*{-1cm}
\subfigure[]{\includegraphics[trim=35mm 15mm 40mm 0mm,clip,scale=0.60]{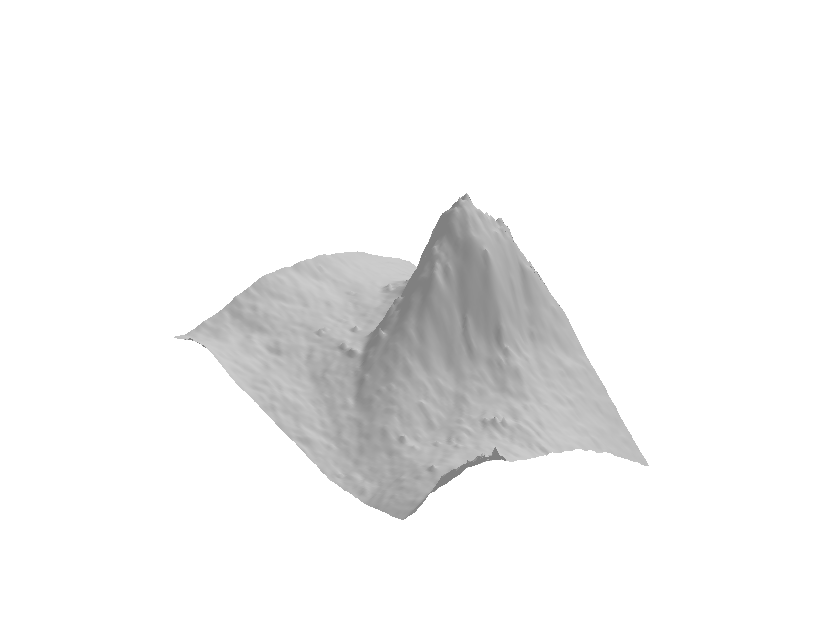}}
\hspace*{-.75cm}
\subfigure[]{\includegraphics[trim=40mm 15mm 30mm 0mm,clip,scale=0.60]{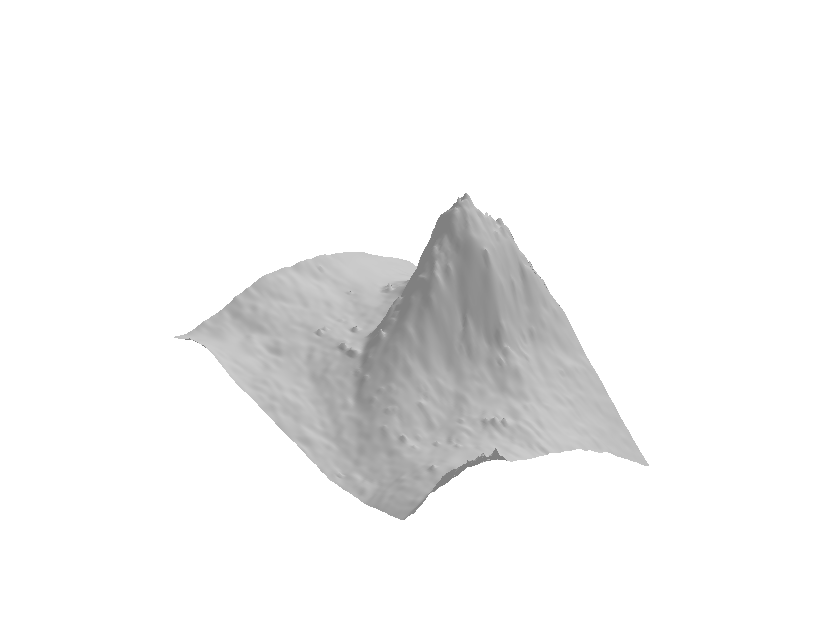}}\\
\hspace*{-.75cm}
\subfigure[]{\includegraphics[trim=0mm 0mm 0mm 0mm,clip,scale=0.40]{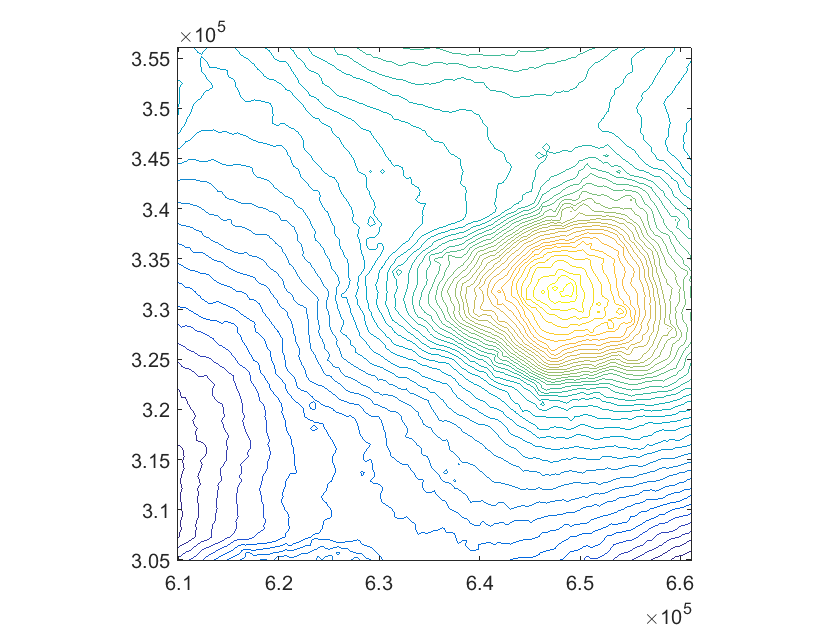}}
\hspace*{-.75cm}
\subfigure[]{\includegraphics[trim=0mm 0mm 0mm 0mm,clip,scale=0.40]{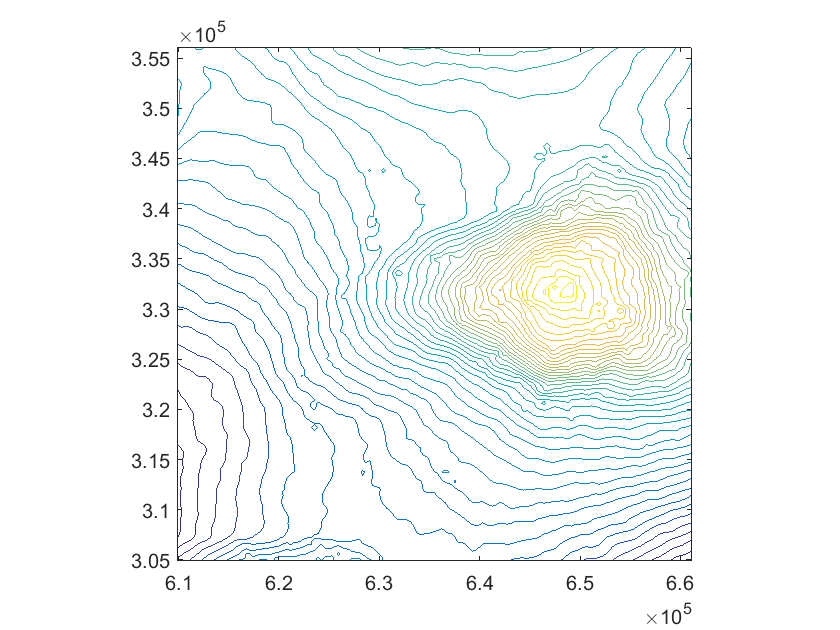}}
\caption{The original data surface of Example \ref{ex4} (a) and the hierarchical spline approximation (b). The contour plots of the original data (c) and of our approximation (d) are also shown.}
 \label{ha2}
%\end{center}
\end{figure}

\end{exm}

The adaptive nature of Algorithm~1 is outlined in Table~\ref{tab:perc} where the information concerning the degrees of local polynomials for Examples~1--4 based on biquadratic splines is reported.

\begin{table}[H]
\centerline{ 
\begin{footnotesize}\begin{tabular}{SSSSS} \toprule
 & {NDOF} & {degree 0} & {degree 1} & {degree 2} \\ \midrule
{example 1}  & {$23092$} & 12.186\% & 24.043\% & 63.771\% \\
{example 2}  &   {$2736$} & 18.469\% & 22.348\% & 59.183\% \\
{example 3}  & {$18218$} & 23.234\% & 29.793\% & 46.973\% \\
{example 4}  & {$59706$} & 47.459\% & 35.413\% & 17.127\% \\
\bottomrule
\end{tabular}\end{footnotesize}
}
\caption{Percentage of polynomials of degrees 0, 1, 2 computed with Algorithm~1 for Examples 1-4.}\label{tab:perc}
\end{table}

In these experiments, biquadratic splines proved to be the best compromise between the need to keep the approximation local and the possibility of considering polynomials of a suitable (high) degree.

\begin{exm}\label{ex5}
In the last example, the set $F$ defined in \eqref{scatset} is obtained by sampling the peak function
$f(x,y) = {2}/\left({3\,\exp((10\, x - 3)^2+(10\, y + 3)^2)}\right)
$ for $(x,y)\in[-1,1]^2$ on the set of $16000$ scattered data locations shown in Figure~\ref{fig:p1} (left). When considering $\sigma=10^{-6}$, $\epsilon = 2\cdot 10^{-3}$, $M_{max} = 7$, and an initial uniform mesh with $15\times 15$ elements, with ${\mathbf d} = (2,2)$ the scheme does not meet the given tolerance $\epsilon$.
%, since from level 5 we have a deterioration in the convergence behaviour. 
This problem is resolved by considering ${\mathbf d} = (4,4)$ which allows us to obtain the final approximation with $2390$ degrees of freedom distributed on four levels. Polynomials of degree 4 are always selected in this case. The corresponding hierarchical mesh is shown in Figure~\ref{fig:p1} (right). Clearly, the situation can remarkably vary according to the initialization of the input parameters.

\begin{figure}[h!]%\begin{center}
\centering{\hspace*{-1.9cm}
\subfigure[]{\includegraphics[trim=0mm 0mm 0mm 0mm,clip,scale=0.50]{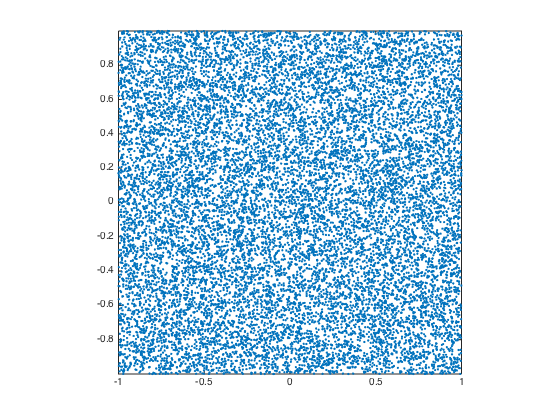}}
\subfigure[]{\raisebox{7mm}{{\includegraphics[trim=0mm 0mm 0mm 0mm,clip,scale=0.50]{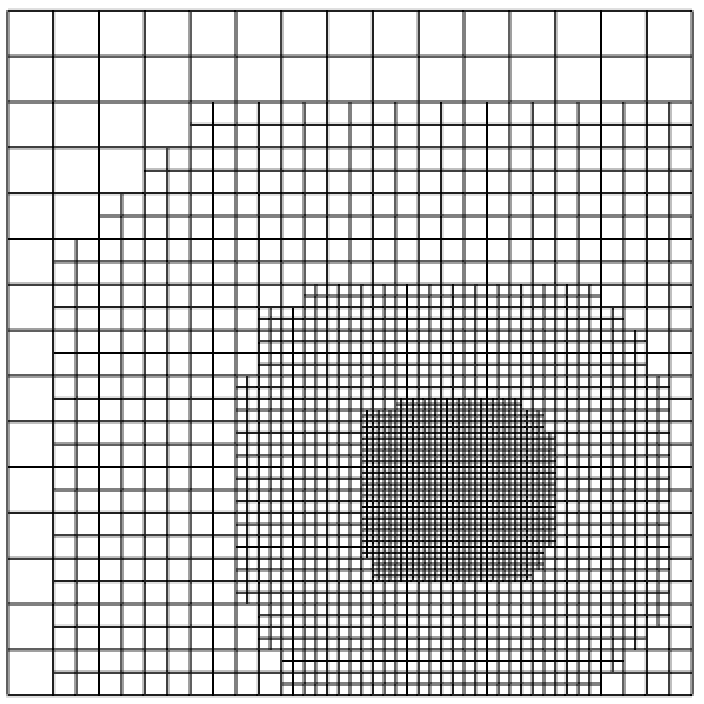}}
}}}
\caption{The data point (a) and the hierarchical mesh (b) generated by our approximation algorithm in Example \ref{ex5} for $\mathbf{d}=(4,4)$.}
 \label{fig:p1}
%\end{center}
\end{figure}

\end{exm}

\begin{figure}[h!]%\begin{center}
\centerline{
\subfigure[]{\includegraphics[trim=0mm 0mm 0mm 0mm,clip,scale=0.50]{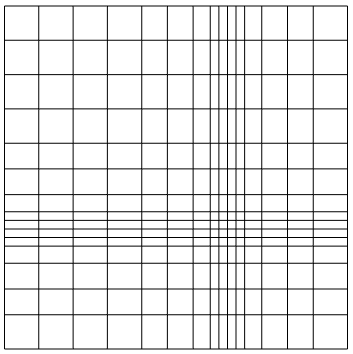}}\hspace*{1.6cm}
\subfigure[]{\includegraphics[trim=0mm 0mm 0mm 0mm,clip,scale=0.50]{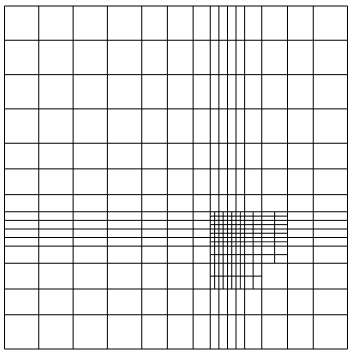}}
}
\caption{Initial non-uniform mesh (a) and hierarchical mesh on two levels (b) in Remark~\ref{rmk:nonuni}.}
 \label{fig:nonuni}
%\end{center}
\end{figure}

\begin{remark}\label{rmk:nonuni}Both the theoretical framework of (truncated) hierarchical B-splines and the corresponding algorithms here presented are not restricted to uniform knot configurations but also cover the case of arbitrary knot distributions. Since the adaptive nature of the mesh is already realized with the hierarchical construction, the common practice working with (T)HB-splines relies on starting with a uniform tensor-product grid, for then successively identifying the different levels of resolution via dyadic refinement. Nevertheless, it is possible to start with a non-uniform initial grid and refine every marked cell by performing an arbitrary splitting into subcells. Obviously, the resulting implementation can become unnecessarily complicated, also due to the arbitrary choice of the splitting. As a compromise, we consider a test where we start with a non-uniform mesh for then performing dyadic refinement. Figure~\ref{fig:nonuni} (b) shows the mesh obtained with $\mathbf{d}=(2,2)$ for the data set considered in Example~5 starting with the non-uniform mesh with $15\times 15$ elements reported on Figure~\ref{fig:nonuni} (a). In this case, by consdering $\sigma = 10^{-6}$ and $\epsilon = 5\cdot 10^{-2}$, the tolerance $\epsilon$ is met in two levels with $370$ degrees of freedom. Repeating the same experiment starting with a uniform mesh with the same number of elements, we need three levels and $553$ degrees of freedom.
\end{remark}

\section{Conclusion}\label{sec:end}

An adaptive scattered data fitting scheme based on hierarchical spline spaces has been introduced by defining a quasi-interpolant whose coefficients are obtained from local polynomial approximations of the data. To this aim, local least squares  polynomial approximations of variable degree are suitably combined with hierarchical quasi-interpolation. By exploiting the characterization of THB-spline constructions, the coefficients obtained as solution of the local problems can be used to define the global spline approximation without the need of additional computations. The resulting scheme has been validated on several scattered data sets of different nature including structured data, data with highly varying density, data with holes, and also configurations along contour lines. The results confirm that high quality approximations can be directly computed via \emph{compact} hierarchical representations which need a reduced number of degrees of freedom compared to previous two-stage methods.

%%%%%%%%%%%%%%%%%%%%%%%%%%%%%%%%%%%%%%%%%%%%%%%%%%%%%%%%%%%%%%%%%%%%
\section*{Acknowledgements}
The support by MIUR ``Futuro in Ricerca" programme through the project DREAMS (RBFR13FBI3) and by the Istituto Nazionale di Alta Matematica (INdAM) through Gruppo Nazionale per il Calcolo Scientifico (GNCS)  are gratefully acknowledged. 

\bibliographystyle{abbrv}
\bibliography{biblio}

\begin{thebibliography}{10}

\bibitem{bracco2012}
G.~Allasia and C.~Bracco.
\newblock {Multivariate Hermite-Birkhoff interpolation by a class of cardinal
  basis functions}.
\newblock {\em Appl. Math. Comput.}, 218:9248--9260, 2012.

\bibitem{bracco2016}
C.~Bracco, C.~Giannelli, F.~Mazzia, and A.~Sestini.
\newblock {Bivariate hierarchical Hermite spline quasi--interpolation}.
\newblock {\em \BIT}, 56:1165--1188, 2016.

\bibitem{davydov2002}
O.~Davydov.
\newblock On the approximation power of local least squares polynomials.
\newblock In J.~Levesley, I.~J. Anderson, and J.~C. Mason, editors, {\em
  {Algorithms for Approximation IV}}, pages 346--353. University of
  Huddersfield, UK, 2002.

\bibitem{davydov2005}
O.~Davydov, R.~Morandi, and A.~Sestini.
\newblock {Local RBF approximation for scattered data fitting with bivariate
  splines}.
\newblock In M.~De~Bruin, D.~H. Mache, and J.~Szabados, editors, {\em Trends
  and Applications in Constructive Approximation}, volume 151 of {\em
  International Series of Numerical Mathematics}, pages 91--102. Birkh\"auser /
  Basel, 2005.

\bibitem{davydov2006}
O.~Davydov, R.~Morandi, and A.~Sestini.
\newblock Local hybrid approximation for scattered data fitting with bivariate
  splines.
\newblock {\em \CAGD}, 23:703--721, 2006.

\bibitem{davydov2014}
O.~Davydov, J.~Prasiswa, and U.~Reif.
\newblock Two-stage approximation methods with extended b-splines.
\newblock {\em \MC}, 83:809--833, 2014.

\bibitem{davydov2004}
O.~Davydov and F.~Zeilfelder.
\newblock Scattered data fitting by direct extension of local polynomials to
  bivariate splines.
\newblock {\em \ACM}, 21:223--271, 2004.

\bibitem{fasshauer2007}
G.~F. Fasshauer.
\newblock {\em Meshfree Approximation Methods with MATLAB}.
\newblock World Scientific Publishing Co., Inc., River Edge, NJ, USA, 2007.

\bibitem{forsey1988}
D.~R. Forsey and R.~H. Bartels.
\newblock Hierarchical {B}-spline refinement.
\newblock {\em \CG}, 22:205--212, 1988.

\bibitem{forsey1995}
D.~R. Forsey and R.~H. Bartels.
\newblock Surface fitting with hierarchical splines.
\newblock {\em \ACMTG}, 14:134--161, 1995.

\bibitem{giannelli2012}
C.~Giannelli, B.~J\"uttler, and H.~Speleers.
\newblock {THB}-splines: The truncated basis for hierarchical splines.
\newblock {\em \CAGD}, 29:485--498, 2012.

\bibitem{giannelli2014}
C.~Giannelli, B.~J\"uttler, and H.~Speleers.
\newblock Strongly stable bases for adaptively refined multilevel spline
  spaces.
\newblock {\em \ACM}, 40:459--490, 2014.

\bibitem{greiner1997}
G.~Greiner and K.~Hormann.
\newblock Interpolating and approximating scattered {3D}-data with hierarchical
  tensor product {B}-splines.
\newblock In A.~L. M{\'e}haut{\'e}, C.~Rabut, and L.~L. Schumaker, editors,
  {\em Surface Fitting and Multiresolution Methods}, pages 163--172. Vanderbilt
  University Press, Nashville, TN, 1997.

\bibitem{hjelle2005}
{\O}.~Hjelle and M.~D{\ae}hlen.
\newblock Multilevel least squares approximation of scattered data over binary
  triangulations.
\newblock {\em Computing and Visualization in Science}, 8(2):83--91, 2005.

\bibitem{kiss2014b}
G.~Kiss, C.~Giannelli, U.~Zore, B.~J\"uttler, D.~Gro{\ss}mann, and J.~Barner.
\newblock {Adaptive CAD model (re-)construction with {THB}-splines}.
\newblock {\em \GM}, 76:273--288, 2014.

\bibitem{kraft1997}
R.~Kraft.
\newblock Adaptive and linearly independent multilevel {B}-splines.
\newblock In A.~Le~M{\'e}haut{\'e}, C.~Rabut, and L.~L. Schumaker, editors,
  {\em Surface Fitting and Multiresolution Methods}, pages 209--218. Vanderbilt
  University Press, Nashville, 1997.

\bibitem{kraft1998}
R.~Kraft.
\newblock {\em Adaptive und linear unabh\"angige Multilevel {B}-Splines und
  ihre Anwendungen}.
\newblock PhD thesis, Universit\"at Stuttgart, 1998.

\bibitem{lai2009}
M.-J. Lai and L.~L. Schumaker.
\newblock A domain decomposition method for computing bivariate spline fits of
  scattered data.
\newblock {\em \SIAMJNA}, 47:911--928, 2009.

\bibitem{lee2005}
B.-G. Lee, J.-J. Lee, and K.-R. Kwon.
\newblock {Quasi-interpolants based multilevel B-spline surface reconstruction
  from scattered data}.
\newblock In O.~Gervasi, M.~L. Gavrilova, V.~Kumar, A.~Lagan{\`a}, H.~P. Lee,
  Y.~Mun, D.~Taniar, and C.~J.~K. Tan, editors, {\em Computational Science and
  Its Applications -- ICCSA 2005: International Conference, Singapore, May
  9-12, 2005, Proceedings, Part III}, pages 1209--1218. Springer Berlin
  Heidelberg, 2005.

\bibitem{liu2012}
S.~Liu and C.~L. Wang.
\newblock Quasi-interpolation for surface reconstruction from scattered data
  with radial basis functions.
\newblock {\em \CAGD}, 29:435--447, 2012.

\bibitem{nurnberger2005}
G.~N\"urnberger, V.~Rayevskaya, L.~L. Schumaker, and F.~Zeilfelder.
\newblock Local lagrange interpolation with bivariate splines of arbitrary
  smoothness.
\newblock {\em \CA}, 23:33--59, 2005.

\bibitem{ohtake2006}
Y.~Ohtake, A.~Belyaev, and H.-P. Seidel.
\newblock Sparse surface reconstruction with adaptive partition of unity and
  radial basis functions.
\newblock {\em \GM}, 68:15--24, 2006.

\bibitem{rabut2005}
C.~Rabut.
\newblock Locally tensor product functions.
\newblock {\em \NA}, 39:329--348, 2005.

\bibitem{skytt2015}
V.~Skytt, O.~Barrowclough, and T.~Dokken.
\newblock Locally refined spline surfaces for representation of terrain data.
\newblock {\em \CG}, 49:58--68, 2015.

\bibitem{speleers2015b}
H.~Speleers.
\newblock {A family of smooth quasi-interpolants defined over Powell-Sabin
  triangulations}.
\newblock {\em \CA}, 41:297--324, 2015.

\bibitem{speleers2017}
H.~Speleers.
\newblock Hierarchical spline spaces: quasi-interpolants and local
  approximation estimates.
\newblock {\em \ACM}, 43:235--255, 2017.

\bibitem{speleers2016}
H.~Speleers and C.~Manni.
\newblock Effortless quasi-interpolation in hierarchical spaces.
\newblock {\em \NM}, 132:155--184, 2016.

\bibitem{usgeo}
U.~S.~G. Survey.
\newblock \url{https://www.usgs.gov/},
  \url{http://dds.cr.usgs.gov/pub/data/nationalatlas/el_usa_hawaii.bil_nt00924.tar.gz}.

\bibitem{vuong2011}
A.-V. Vuong, C.~Giannelli, B.~J\"uttler, and B.~Simeon.
\newblock A hierarchical approach to adaptive local refinement in isogeometric
  analysis.
\newblock {\em \CMAME}, 200:3554--3567, 2011.

\bibitem{wendland2004}
H.~Wendland.
\newblock {\em Scattered Data Approximation}.
\newblock Cambridge Monographs on Applied and Computational Mathematics.
  Cambridge University Press, 2004.

\end{thebibliography}

\end{document}